\documentclass[12pt,psfig]{article}
\usepackage{graphicx,epsfig}
\usepackage{amsmath}
\usepackage{mathrsfs}
\usepackage{amssymb}
\def\bbb{\begin{eqnarray*}}

\def\eee{\end{eqnarray*}}

\begin{document}

\baselineskip=17pt
\begin{center}

\vspace{-0.6in} {\large \bf Topological conjugacy between induced non-autonomous set-valued systems and subshifts of finite type}\\ [0.2in]

Hua Shao$^\dag$, Guanrong Chen$^\dag$, Yuming Shi$^\ddag$

\vspace{0.15in}$^\dag$ Department of Electronic Engineering, City University of Hong Kong,\\
Hong Kong SAR, P.~R. China\\

\vspace{0.15in}$^\ddag$ Department of Mathematics, Shandong University \\
 Jinan, Shandong 250100, P.~R. China\\
 
\end{center}

\footnote{$^\ddag$ The corresponding author: ymshi@sdu.edu.cn}

{\bf Abstract.} This paper establishes topological (equi-)semiconjugacy and (equi-)conjugacy between induced
non-autonomous set-valued systems and subshifts of finite type. First, some necessary and sufficient conditions are given
for a non-autonomous discrete system to be topologically semiconjugate or conjugate to a subshift of finite type. Further, several
sufficient conditions for it to be topologically equi-semiconjugate or equi-conjugate to a subshift of finite type are obtained.
Consequently, estimations of topological entropy and several criteria of Li-Yorke chaos and distributional chaos in a sequence
are derived. Second, the relationships of several related dynamical behaviors between the non-autonomous discrete system and
its induced set-valued system are investigated. Based on these results, the paper furthermore establishes the topological
(equi-)semiconjugacy and (equi-)conjugacy between induced set-valued systems and subshifts of finite type. Consequently,
estimations of the topological entropy for the induced set-valued system are obtained, and several criteria of Li-Yorke chaos and
distributional chaos in a sequence are established. Some of these results not only extend the existing related results for
autonomous discrete systems to non-autonomous discrete systems, but also relax the assumptions of the counterparts in the literature.
Two examples are finally provided for illustration.
\medskip

\medskip

{2010 {\bf \it Mathematics Subject Classification}}: 37B55, 54C60, 37B10.

\bigskip

\noindent{\bf 1. Introduction}\medskip

Symbolic dynamical systems play a significant role in the study of chaos theory of dynamical systems since
they appear to be simple but have quite rich and complex dynamical behaviors. They were first utilized
by Hadamard to study geodesics on surfaces with negative curvatures \cite{hadamard}. Then, they were applied
by Birkhoff to study dynamical systems \cite{birkhoff}. In the 1960s, complex dynamical behaviors of Smale
horseshoe were depicted by its topological conjugacy to a two-sided symbolic  system \cite{smale}. Two
well-known maps, the H\'enon map and the logistic map, were shown to be chaotic for some parameters also
by their topological conjugacy to the symbolic dynamical systems \cite{Devaney79, Devaney89, Robinson}.
In fact, it is a very useful method for studying the complexity of dynamical systems by establishing topological
semiconjugacy or conjugacy to symbolic dynamical systems (see \cite{block,Kennedy,Shi04,Shi06,Shi09,shicsf,zhang},
and references therein).

It is worth mentioning that Block and Coppel in 1992 introduced the concept of ``turbulence" for continuous interval
maps, and proved that a strictly turbulent map is topologically semi-conjugate to the one-sided symbolic system on two
symbols in a compact invariant set \cite{block}. In 2004, we proved that a strictly turbulent map, which satisfies
an expanding condition in distance, is topologically conjugate to the one-sided symbolic system on two symbols \cite{Shi04}.
In 2006, we changed the term  ``turbulence" to ``coupled-expansion" in order to avoid possible confusion with the
turbulence concept in fluid mechanics \cite{Shi06}. Later, we extended this notion to coupled-expansion for a transition
matrix and showed that a map, which is strictly coupled-expanding for a transition matrix, is topological conjugate to
a subshift of finite type under certain conditions \cite{shicsf}. Thereafter, chaos induced
by coupled-expansion for a transition matrix has attracted interest and attention from researchers in the field \cite{ju16,Kim,Kulczycki,Shao15,Shi09,Zhangli,zhangijbc,zhang}.

In the study of complex dynamics, sometimes it is not enough to know only the tracjectory of a single point, but it needs
to know the motion of a collection of points. For instance, in the study of collective behaviors in biology, one needs to
know the massive migration of birds or mammals. Inspired by this natural phenomenon, Rom\'an-Flores studied the
relationship between transitivity of a continuous map and its induced set-valued system \cite{Rom}.
Following his work, many scholars studied the relationships between individual chaos and collective chaos
\cite{Fed,Gu06,ju17,Liao,Liu,Ma,san,wang,wu}. In particular, Wang and Wei proved that, if a
continuous map is strictly coupled-expanding, then its induced set-valued discrete system is
topologically semi-conjugate to a full shift \cite{wang}. Recently, Ju et al. generalized
this result and gave some sufficient conditions under which the induced set-valued system is
topologically (semi-)conjugate to a subshift of finite type \cite{ju17}.

Motivated by the above research, we are interested in studying the topological (equi-)semiconjugacy and (equi-)conjugacy
between the induced non-autonomous set-valued system and a subshift of finite type. We achieve this goal in two steps.
First, we establish the topological (equi-)semiconjugacy and (equi-)conjugacy between a non-autonomous
discrete system and a subshift of finite type. Second, we investigate the relationships of several related
dynamical behaviors between the non-autonomous discrete system and its induced set-valued system. Based
on these results, we can establish the topological (equi-)semiconjugacy and (equi-)conjugacy between the
induced non-autonomous set-valued system and a subshift of finite type. Since the complex dynamical behaviors of
subshifts of finite type are well understood (see, for example, \cite{Zhou}), and the
relationships of Li-Yorke chaos and topological entropy of two topological equi-(semi)conjugate non-autonomous discrete
systems have been given in \cite{Shi09} and \cite{kolyada,Shao16}, respectively, to that end several criteria of
Li-Yorke chaos and distributional chaos in a sequence and estimations of the topological entropy for the induced
set-valued system are obtained in this paper.

The rest of the present paper is organized as follows. Section 2 presents some basic concepts and useful lemmas.
In Section 3, the topological (equi-)semiconjugacy and (equi-)conjugacy  between a non-autonomous discrete
system and a subshift of finite type are studied, and consequently estimations of topological entropy and some
criteria of Li-Yorke chaos and distributional chaos in a sequence for the non-autonomous discrete
system are established. In Section 4, the relationships of several related dynamical behaviors between
the non-autonomous discrete system and its induced set-valued system are discussed, and the topological
(equi-)semiconjugacy and (equi-)conjugacy between the induced non-autonomous set-valued system and a subshift
of finite type are investigated. By applying the above-obtained results, estimations of the topological entropy for
the induced non-autonomous set-valued system are obtained, and several criteria of Li-Yorke chaos and
distributional chaos in a sequence are established. Two examples are finally provided in Section 5 for illustration.

\bigskip

\noindent{\bf 2. Preliminaries}\medskip

This section is divided into three parts for convenience of discussion. In Section 2.1, several basic concepts and related lemmas
are introduced. Then, the notion of induced non-autonomous set-valued system is introduced in Section 2.2. In Section 2.3, the concepts
of subshifts of finite type and weak coupled-expansion for transition matrices are reviewed.\medskip

\noindent{\bf 2.1. Some basic concepts and related lemmas}\medskip

Consider the following non-autonomous discrete system:
\vspace{-0.2cm}$$ x_{n+1} = f_n(x_n), \;\;n \geq0,                                                                                   \eqno(2.1)\vspace{-0.2cm}$$
where $(X,d)$ is a metric space and $f_n: X\to X$ is a map, $n\geq0$. Denote
$f_{0,\infty}:=\{f_n\}_{n=0}^{\infty}$, $f^n_i:=f_{i+n-1}\circ\cdots\circ f_i$, and $f^{-n}_i:=(f_{i}^n)^{-1}$,
$i\geq0$, $n\geq1$. Let $x\in X$ and $A,B$ be nonempty subsets of $X$. The boundary of $A$ is denoted by $\partial A$;
the diameter of $A$ is denoted by $d(A)$; the distance between $x$ and $A$ is denoted by $d(x,A):=\inf\{d(x,a): a\in A\}$;
and the distance between $A$ and $B$ is denoted by $d(A,B):=\inf\{d(a,b): a\in A,b\in B\}$. The set of all nonnegative integers
and positive integers are denoted by $\mathbf{N}$ and $\mathbf{Z}^{+}$, respectively.
\medskip

\noindent{\bf Definition 2.1} (\cite{Shi09}, Definition 2.7).
System (2.1) is said to be Li-Yorke $\delta$-chaotic for some $\delta>0$ if it has an uncountable Li-Yorke
$\delta$-scrambled set $S$ in $X$; that is, for any $x, y\in S\subset X$,
\vspace{-0.2cm}$$\liminf_{n\to\infty}d(f_{0}^{n}(x),f_{0}^{n}(y))=0\;{\rm and}\;\limsup_{n\to\infty}d(f_{0}^{n}(x),f_{0}^{n}(y))>\delta.\vspace{-0.2cm}$$
Further, it is said to be chaotic in the strong sense of Li-Yorke if all the orbits starting from the points in $S$ are bounded.\medskip

\noindent{\bf Definition 2.2} (\cite{shaocsf}, Definitions 2.1 and 2.2).
System (2.1) is said to be distributionally chaotic if it has an uncountable distributionally scrambled set $D$ in $X$; that is, for any $x, y\in D\subset X$,
\begin{itemize}\vspace{-0.2cm}
\item[{\rm (i)}] $\limsup_{n\to\infty}\frac{1}{n}\sum_{i=0}^{n-1}\chi_{[0,\epsilon)}\big(d(f_{0}^{i}(x),f_{0}^{i}(y))\big)=1$ for any $\epsilon>0$,\vspace{-0.2cm}
\item[{\rm (ii)}] $\liminf_{n\to\infty}\frac{1}{n}\sum_{i=0}^{n-1}\chi_{[0,\delta)}\big(d(f_{0}^{i}(x),f_{0}^{i}(y))\big)=0$ for some $\delta>0$,
\end{itemize}
where $\chi_{[0,\epsilon)}$ is the characteristic function defined on the set $[0,\epsilon)$.\medskip

The definition of topological entropy for system (2.1) is introduced in \cite{kolyada}.
Let $X$ be compact, $Y$ be a nonempty subset of $X$, $\mathcal{A}$ be an open cover of $X$, and $\mathcal{N}(\mathcal{A})$
be the minimal possible cardinality of all subcovers chosen from $\mathcal{A}$. Denote the cover $\{A\cap Y: A\in\mathcal{A}\}$
of the set $Y$ by $\mathcal{A}|_{Y}$ and $\mathcal{A}_{0}^{n}:=\{\bigcap_{j=0}^{n-1}f_{0}^{-j}(A_{j}): A_{j}\in\mathcal{A},\;0\leq j\leq n-1\}$. Then the topological entropy of system (2.1) on $Y$ is defined by
\vspace{-0.2cm}$$h(f_{0,\infty},Y):=\sup_{\mathcal{A}}\left\{\limsup_{n\to\infty}\log\mathcal{N}(\mathcal{A}_{0}^{n}|_{Y})/n\right\}.\vspace{-0.2cm}$$
If $Y=X$, then $h(f_{0,\infty}):= h(f_{0,\infty},X)$ is called the topological entropy of system (2.1) on $X$.\medskip

Now, recall the definitions of topological (equi-)semiconjugacy and (equi-)conjugacy between system (2.1) and
\vspace{-0.2cm}$$u_{n+1}=g_n(u_n),\; n\geq0,                                                                                         \eqno(2.2)\vspace{-0.2cm}$$
defined on a metric space $(Y,e)$, where $g_n: Y\to Y$ is a map, $n\geq0$.\medskip

Let $\{\Lambda_n\}_{n=0}^{\infty}$ and $\{E_n\}_{n=0}^{\infty}$ be two sequences of subsets of $X$ and $Y$,
respectively, and $h_n: \Lambda_n\to E_n$ be a map, $n\geq0$. The sequence of maps $\{h_n\}_{n=0}^{\infty}$
is called equi-continuous in $\{\Lambda_n\}_{n=0}^{\infty}$ if for any $\epsilon>0$ there exists $\delta>0$ such that
$e(h_n(x), h_n(y))<\epsilon$ for any $n\geq0$ and any $x,y\in\Lambda_n$ with $d(x,y)<\delta$. In addition,
$\{\Lambda_n\}_{n=0}^{\infty}$ is called invariant under system (2.1) if $f_n(\Lambda_n)\subset\Lambda_{n+1}$
for any $n\geq0$. Then, system (2.1) restricted to $\{\Lambda_n\}_{n=0}^{\infty}$ is an invariant subsystem
of system (2.1) on $\{\Lambda_n\}_{n=0}^{\infty}$ \cite{Shao16}.\medskip

\noindent{\bf Definition 2.3} (\cite{Shi09}, Definition 3.3).
Let $\{\Lambda_n\}_{n=0}^{\infty}$ and $\{E_n\}_{n=0}^{\infty}$ be invariant under systems (2.1) and {\rm(2.2)},
respectively. If for any $n\geq0$, there exists an {\rm(}equi-{\rm)}continuous surjective map $h_n: \Lambda_n\to E_n$
such that $h_{n+1}\circ f_n=g_n\circ h_n$, then the invariant subsystem of system
system (2.1) on $\{\Lambda_n\}_{n=0}^{\infty}$ is said to be topologically $\{h_n\}_{n=0}^{\infty}$
-{\rm(}equi-{\rm)}semiconjugate to the invariant subsystem of system {\rm(2.2)} on $\{E_n\}_{n=0}^{\infty}$.
Further, if $h_n$ is invertible for all $n\geq0$ and $\{h_{n}^{-1}\}_{n=0}^{\infty}$ is also
{\rm(}equi-{\rm)}continuous, then the above two invariant subsystems are said to be topologically
$\{h_n\}_{n=0}^{\infty}$-{\rm(}equi-{\rm)}conjugate.\medskip

The relationships of Li-Yorke chaos and distributional chaos between two topological
equi-conjugate systems are given below, respectively.\medskip

\noindent{\bf Lemma 2.1} (\cite{shaocnsns,Shi09}). {\it Assume that the invariant subsystem of system {\rm(2.1)} on
$\{\Lambda_n\}_{n=0}^{\infty}$ is topologically equi-conjugate to the invariant subsystem of system {\rm(2.2)}
on $\{E_n\}_{n=0}^{\infty}$. Then,
\begin{itemize}\vspace{-0.2cm}
\item[{\rm (i)}] system {\rm(2.1)} has an uncountable Li-Yorke $\delta$-scrambled set in $\Lambda_0$ if and only if system
{\rm(2.2)} has an uncountable Li-Yorke $\gamma$-scrambled set in $E_0$, for some $\delta,\gamma>0$;\vspace{-0.2cm}
\item[{\rm (ii)}] system {\rm(2.1)} has an uncountable distributionally scrambled set in $\Lambda_0$ if and only if system
{\rm(2.2)} has an uncountable distributionally scrambled set in $E_0$.\vspace{-0.2cm}
\end{itemize}}\medskip

The next result shows the relationship of the topological entropy between two topological equi-(semi)conjugate systems.\medskip

\noindent{\bf Lemma 2.2} (\cite{Shao16}, Lemma 2.2). {\it Let $(X,d)$ and  $(Y,e)$ be two compact metric spaces.
If an invariant subsystem of system {\rm(2.1)} on $\{\Lambda_n\}_{n=0}^{\infty}$ is topologically
equi-semiconjugate to an invariant subsystem of system {\rm(2.2)} on $\{E_n\}_{n=0}^{\infty}$,
then $h(f_{0,\infty},\Lambda_0)\geq h(g_{0,\infty},E_0)$. Further, if they are topologically equi-conjugate,
then $h(f_{0,\infty},\Lambda_0)=h(g_{0,\infty},E_0)$.}\medskip

The following result will also be useful in the sequent sections.\medskip

\noindent{\bf Lemma 2.3.} {\it Let $\{E_n\}_{n=0}^{\infty}$ be a sequence of nonempty compact subsets of $X$ satisfying that
$E_{n+1}\subset E_n$ for all $n\geq0$. Then, $\bigcap_{n=0}^{\infty}{E_{n}}$ is a singleton if and only if $\lim_{n\to \infty}d(E_n)=0$.}\medskip

\noindent{\bf Proof.} The sufficiency can be directly derived form Lemma 2.7 in \cite{Shi04}. To show the necessity,
assume otherwise. Then, $\lim_{n\to\infty}d(E_n)=\epsilon>0$. So, there exists $N\in\mathbf{Z}^{+}$ such that $d(E_n)>\epsilon/2$
for all $n\geq N$. Thus, there exist $x_n,y_n\in E_n$ satisfying that
\vspace{-0.2cm}$$d(x_n,y_n)>\frac{\epsilon}{2}, \;n\geq N.                                                                           \eqno(2.3)\vspace{-0.2cm}$$
Fix any $n\geq 0$, and let $m\geq\max\{N,n\}$. Then, $\{x_k\}_{k=m}^{\infty},\{y_k\}_{k=m}^{\infty}
\subset E_m\subset E_n$. Suppose that $x_n\to x$ and $y_n\to y$ as $n\to\infty$. Since $E_n$ is compact, $x,y\in E_n$.
Thus, $x,y\in\bigcap_{n=0}^{\infty}E_n$. By (2.3) one has that $x\neq y$. This is a contradiction to the assumption
that $\bigcap_{n=0}^{\infty}{E_{n}}$ is a singleton. Therefore, $\lim_{n\to \infty}d(E_n)=0$.
This completes the proof.\medskip

\noindent{\bf 2.2. The induced non-autonomous set-valued system}\medskip

Let $\mathcal{K}(X)$ be the class of all nonempty compact subsets of $X$ and the Hausdorff metric on $\mathcal{K}(X)$
be defined by
\vspace{-0.2cm}$$\mathcal{H}_{d}(A, B):=\max\bigg\{\sup_{a\in A}d(a, B), \sup_{b\in B}d(b, A)\bigg\},                                \eqno(2.4)\vspace{-0.2cm}$$
where $A,B\in \mathcal{K}(X)$. It is known that $(\mathcal{K}(X),\mathcal{H}_{d})$
is a compact metric space if and only if $(X, d)$ is a compact metric space \cite{Nadler}.
For convenience, let $\mathcal{H}_{d}(\mathcal{U})$ denote the diameter of a subset $\mathcal{U}$
of $\mathcal{K}(X)$, and let $K$ be a nonempty subset of $X$. Denote
\vspace{-0.2cm}$$\langle K\rangle:=\{A\in\mathcal{K}(X): A\subset K\}.                                                               \eqno(2.5)\vspace{-0.2cm}$$
Let $f_n$ be continuous in $X$, $n\geq0$. Then, system (2.1) induces the following non-autonomous set-valued system:
\vspace{-0.2cm}$$A_{n+1}=\bar{f}_n(A_n),\;\;n\geq0,                                                                                  \eqno(2.6)\vspace{-0.2cm}$$
where $\bar{f}_n$ is defined by $\bar{f}_n(A):=f_n(A)$ for all $A\in\mathcal{K}(X)$. Thus, $f_n(A)\in\mathcal{K}(X)$
and $\bar{f}_n$ is a map from $\mathcal{K}(X)$ to $\mathcal{K}(X)$, $n\geq0$.
Denote $\bar{f}_{0,\infty}:=\{\bar{f}_n\}_{n=0}^{\infty}$.\medskip

\noindent{\bf Lemma 2.4.} {\rm(\cite{shao19}, Lemmas II.3.)} {\it $\{\bar{f}_n\}_{n=0}^{\infty}$ is equi-continuous in $\mathcal{K}(X)$
if and only if $\{f_n\}_{n=0}^{\infty}$ is equi-continuous in $X$. Consequently, $\bar{f}:\mathcal{K}(X)\to\mathcal{K}(X)$ is continuous
 if and only if $f: X\to X$ is continuous.}\medskip

The following two lemmas will also be useful in the sequel.\medskip

\noindent{\bf Lemma 2.5} {\rm(\cite{ju17}, Lemmas 3.1.1 and 3.1.3).} {\it
\begin{itemize}\vspace{-0.2cm}
\item[{\rm (i)}] Let $X$ be compact and $K\in\mathcal{K}(X)$. Then, $\langle K\rangle$ is a nonempty compact subset of $\mathcal{K}(X)$,
where $\langle K\rangle$ is specified in {\rm(2.5)}.\vspace{-0.2cm}
\item[{\rm (ii)}] Let $K_n\in \mathcal{K}(X)$, $n\geq1$.
Then, $\langle\bigcap_{n=1}^{\infty}K_n\rangle=\bigcap_{n=1}^{\infty}\langle K_n\rangle$.\vspace{-0.2cm}
\end{itemize}}\medskip

\noindent{\bf Lemma 2.6.} {\it Let $X$ be compact, $f_n$ be continuous in $X$, $n\geq0$, and $K\in\mathcal{K}(X)$.
Then, $\langle f_{m}^{-i}(K)\rangle\subset\bar{f}_{m}^{-i}(\langle K\rangle)$ for all $m,i\geq0$.}\medskip

\noindent{\bf Proof.} Fix any $m,i\geq0$. For any $K_0\in \langle f_{m}^{-i}(K)\rangle$, one has that $K_0\in\mathcal{K}(X)$ and
$f_{m}^{i}(K_0)\subset K$. Since $f_n$ is continuous in $X$, $n\geq0$, one has $f_{m}^{i}(K_0)\in\mathcal{K}(X)$.
Thus, $f_{m}^{i}(K_0)\in \langle K\rangle$, and so $K_0\in \bar{f}_{m}^{-i}(\langle K\rangle)$.
Hence, $\langle f_{m}^{-i}(K)\rangle\subset\bar{f}_{m}^{-i}(\langle K\rangle)$.
This completes the proof.\medskip

\noindent{\bf 2.3. Subshifts of finite type and weak-coupled-expansion for transition matrices}\medskip

Recall the definitions of subshifts of finite type \cite{Zhou}.
A matrix $A=(a_{ij})_{N\times N}$ ($N\geq2$) is said to be a transition matrix if $a_{ij}=0$ or $1$ for all $i,j$;
$\sum_{j=1}^{N}a_{ij}\geq1$ for all $i$; and $\sum_{i=1}^{N}a_{ij}\geq1$ for all $j$, $1\leq i,j\leq N$.
A transition matrix $A=(a_{ij})_{N\times N}$ is said to be irreducible if, for each pair $1\leq i, j\leq N$,
there exists $k\in\mathbf{Z}^{+}$ such that $a^{(k)}_{ij}>0$, where $a^{(k)}_{ij}$ denotes the $(i, j)$ entry
of matrix $A^{k}$. For a given transition matrix $A=(a_{ij})_{N\times N}$, denote
\vspace{-0.2cm}$$\Sigma_{N}^{+}(A):=\{s=(s_0,s_1,\cdots): 1\leq s_{j}\leq N, \;a_{s_{j}s_{j+1}}=1, \;j\geq0\}.\vspace{-0.2cm}$$
Note that $\Sigma_{N}^{+}(A)$ is a compact metric space with the metric
\vspace{-0.2cm}$$\hat{\rho}(\alpha, \beta):=\sum_{i=0}^{\infty}\hat{d}(a_i,b_i)/2^{i},\;\;\alpha=(a_0,a_1,\cdots),\;
\beta=(b_0,b_1,\cdots)\in\Sigma_{N}^{+}(A),\vspace{-0.2cm}$$
where $\hat{d}(a_i, b_i)=0$ if $a_i=b_i$, and $\hat{d}(a_i, b_i)=1$ if $a_i\neq b_i$, $i\geq0$.
The map $\sigma_A: \Sigma_{N}^{+}(A)\to\Sigma_{N}^{+}(A)$ with $\sigma_A((s_0,s_1,s_2,\cdots)):=(s_1,s_2,\cdots)$
is said to be a subshift of finite type associated with $A$. Its topological entropy is equal to $\log\rho(A)$,
where
\vspace{-0.2cm}$$\rho(A):=\lim_{n\to\infty}\|A^{n}\|^{\frac{1}{n}},\;\;\|A\|=\sum_{1\leq i,j\leq N}a_{ij}.                           \eqno(2.7)\vspace{-0.2cm}$$
It is known that $(\Sigma_{N}^{+}(A),\sigma_A)$ is Li-Yorke chaotic if and only if it is distributionally chaotic \cite{wanghy}.\medskip

Next, the definition of weak coupled-expansion for a transition matrix is introduced.\medskip

\noindent{\bf Definition 2.4.}
Let $A=(a_{ij})_{N\times N}$ be a transition matrix. If there exists a sequence $\{V_{i,n}\}_{n=0}^{\infty}$ of
nonempty subsets of $X$ with $V_{i,n}\cap V_{j,n}=\partial V_{i,n} \cap\partial V_{j,n}$ {\rm(}$d(V_{i,n},V_{j,n})>0${\rm)}
for all $1\leq i\neq j\leq N$ and $n\geq0$ such that
\vspace{-0.2cm}$$f_n(V_{i,n})\supset\bigcup_{a_{ij}=1}V_{j,n+1},\;1\leq i\leq N, \;n\geq0,\vspace{-0.2cm}$$
then system (2.1) is said to be {\rm(}strictly{\rm)} weakly $A$-coupled
-expanding in $\{V_{i,n}\}_{n=0}^{\infty}$, $1\leq i\leq N$. In the special case that $V_{i,n}=V_i$, $1\leq i\leq N$,
$n\geq0$, it is said to be {\rm(}strictly{\rm)} $A$-coupled-expanding in $V_i$, $1\leq i\leq N$.\medskip

\noindent{\bf Remark 2.1.} Definition 2.4 is a slight revision of that in \cite{Zhangli}.\medskip

Denote
\vspace{-0.2cm}$$
U_i:=\{\alpha=(a_0,a_1,\cdots)\in\Sigma_{N}^{+}(A): a_0=i\},\;1\leq i\leq N.
\eqno(2.8)\vspace{-0.2cm}$$
Then, $U_i$, $1\leq i\leq N$, are disjoint and nonempty compact subsets of $\Sigma_{N}^{+}(A)$ and satisfy that
\vspace{-0.2cm}$$
\bigcup_{i=1}^{N}U_i=\Sigma_{N}^{+}(A).
\eqno(2.9)\vspace{-0.2cm}$$

\noindent{\bf Lemma 2.7} (\cite{shicsf}, Theorem 3.1). {\it $\sigma_{A}$ is strictly $A$-coupled-expanding in $U_i$, $1\leq i\leq N$.}

\bigskip

\noindent{\bf 3. Topological (semi)conjugacy and equi-(semi)conjugacy between non-autonomous discrete systems and subshifts of finite type}\medskip

Now, the topological (semi)conjugacy and equi-(semi)conjugacy between system (2.1) and a subshift of finite type
are studied in Sections 3.1 and 3.2, respectively.\medskip

\noindent{\bf 3.1. Topological semiconjugacy and conjugacy} \medskip

\noindent{\bf Lemma 3.1.} {\it Let $A=(a_{ij})_{N\times N}$ be a transition matrix and $\{V_{i,n}\}_{n=0}^{\infty}$ be a sequence of
nonempty subsets of $X$, $1\leq i\leq N$.
\begin{itemize}
\item[{\rm (i)}] If $f_n$ is continuous in $E$, $n\geq0$, and $V_{i,n}$ is a compact subset of $X$,
$1\leq i\leq N$, $n\geq0$, then $V_{\alpha}^{m,n}$ is a compact subset of $X$ and satisfies that $V_{\alpha}^{m+1,n}\subset V_{\alpha}^{m,n}\subset\bigcup_{i=1}^{N}V_{i,n}$ for all $m,n\geq0$ and all $\alpha=(a_0,a_1,\cdots)\in\Sigma_{N}^{+}(A)$, where
\vspace{-0.2cm}$$
E:=\bigcup_{n=0}^{\infty}\bigcup_{i=1}^{N}V_{i,n},\;\;V_{\alpha}^{m,n}:=\bigcap_{k=0}^{m}f_{n}^{-k}(V_{a_{k},n+k}).
\eqno(3.1)\vspace{-0.2cm}$$
\item[{\rm (ii)}] If system {\rm(2.1)} is weakly $A$-coupled-expanding in $\{V_{i,n}\}_{n=0}^{\infty}$,
$1\leq i\leq N$, then

${\bf(H_1)}$ $V_{\alpha}^{m,n}\neq\emptyset$ for all $m,n\geq0$ and all $\alpha\in\Sigma_{N}^{+}(A)$.

\item[{\rm (iii)}] If $V_i$, $1\leq i\leq N$, are bounded subsets of $X$ with $V_{i,n}\subset V_i$, $n\geq0$,
and there exists $\lambda>1$ such that $d(f_{n}(x),f_{n}(y))\geq \lambda d(x,y)$ for all $x,y\in V_{i,n}$, $1\leq i\leq N$, $n\geq0$, then

${\bf(H_2)}$ $d(V_{\alpha}^{m,n})$ uniformly converges to $0$ with respect to $n\geq0$ as $m\to\infty$ for any $\alpha\in\Sigma_{N}^{+}(A)$.
\item[{\rm (iv)}] If assumptions {\rm (i)}-{\rm (iii)} hold, then $\bigcap_{m=0}^{\infty}V_{\alpha}^{m,n}$ is a singleton
for any $n\geq0$ and any $\alpha\in\Sigma_{N}^{+}(A)$.
\end{itemize}}\medskip

\noindent{\bf Proof.} It can be easily verified that assertions (i) and (ii) hold. Assertion (iv) is a direct consequence of
assertions (i)-(iii) and Lemma 2.3. Thus, it suffices to show assertion (iii). Fix any $\alpha=(a_0,a_1,\cdots)\in\Sigma_{N}^{+}(A)$.
For any $m,n\geq0$ and $x,y\in V_{\alpha}^{m,n}$, $f_{n}^{i}(x),f_{n}^{i}(y)\in V_{a_{i},n+i}$,
$0\leq i\leq m$. This, together with the assumption (iii), yields that
\vspace{-0.2cm}$$d(f_{n}^{m}(x),f_{n}^{m}(y))\geq\lambda^{m}d(x,y),\vspace{-0.2cm}$$
and thus
\vspace{-0.2cm}$$d(x,y)\leq\lambda^{-m}d(f_{n}^{m}(x),f_{n}^{m}(y))\leq\lambda^{-m}\max_{1\leq i\leq N}d(V_i).\vspace{-0.2cm}$$
Consequently,
\vspace{-0.2cm}$$d(V_{\alpha}^{m,n})\leq\lambda^{-m}\max_{1\leq i\leq N}d(V_i).\vspace{-0.2cm}$$
Therefore, $d(V_{\alpha}^{m,n})$ uniformly converges to $0$ with respect to $n\geq0$ as $m\to\infty$.
The proof is complete.\medskip

Next, a necessary and sufficient condition is derived to ensure system (2.1) to have an invariant subsystem
that is topologically semiconjugate to a subshift of finite type. \medskip

\noindent{\bf Theorem 3.1.} {\it Let $A=(a_{ij})_{N\times N}$ be a transition matrix. Then, there exists a sequence
$\{V_{i,n}\}_{n=0}^{\infty}$ of nonempty compact subsets of $X$ with $V_{i,n}\cap V_{j,n}=\emptyset$ for all
$1\leq i\neq j\leq N$ and $n\geq0$ such that $f_n$ is continuous in $E$, $n\geq0$, and satisfies assumption
${\bf(H_1)}$, where $E$ is specified in {\rm(3.1)}, if and only if, for any $n\geq0$, there exists a nonempty compact
subset $\Lambda_{n}\subset X$ with $f_{n}(\Lambda_{n})\subset\Lambda_{n+1}$ such that $f_n$ is continuous in
$\bigcup_{n=0}^{\infty}\Lambda_{n}$ and the invariant subsystem of system {\rm (2.1)} on $\{\Lambda_{n}\}_{n=0}^{\infty}$
is topologically semiconjugate to $(\Sigma_{N}^{+}(A),\sigma_A)$.}\medskip

\noindent{\bf Proof.} {\bf Necessity}. Fix any $n\geq0$. Then, $V_{\alpha}^{m,n}$ is a nonempty compact subset of $X$ and
satisfies that $V_{\alpha}^{m+1,n}\subset V_{\alpha}^{m,n}$, $m\geq0$, $\alpha\in\Sigma_{N}^{+}(A)$,
by (i) of Lemma 3.1. Thus, $\bigcap_{m=0}^{\infty}V_{\alpha}^{m,n}\neq\emptyset$ for all $\alpha\in\Sigma_{N}^{+}(A)$. Denote
\vspace{-0.2cm}$$
\Lambda_{n}:=\bigcup_{\alpha\in\Sigma_{N}^{+}(A)}\bigcap_{m=0}^{\infty}V_{\alpha}^{m,n}.                                             \eqno(3.2)\vspace{-0.2cm}$$
Clearly, $\Lambda_{n}\neq\emptyset$ and $\Lambda_{n}\subset\bigcup_{i=1}^{N}V_{i,n}$. Then, $\bigcup_{n=0}^{\infty}\Lambda_{n}\subset E$,
and thus $f_n$ is continuous in $\bigcup_{n=0}^{\infty}\Lambda_{n}$ since $f_n$ is continuous in $E$.

Suppose that $\{x_k\}_{k=1}^{\infty}\subset\Lambda_{n}$ is a convergent sequence with $x_k\to x$ as $k\to\infty$. Then, for any $k\geq1$,
there exists $\alpha_k\in\Sigma_{N}^{+}(A)$ such that $x_k\in\bigcap_{m=0}^{\infty}V_{\alpha_k}^{m,n}$. Since $\Sigma_{N}^{+}(A)$
is compact, without loss of generality, suppose that $\alpha_k\to\alpha$ as $k\to\infty$. Then, for any $m\geq0$, there exists
$N_m\in\mathbf{Z}^{+}$ such that $V_{\alpha_k}^{m,n}=V_{\alpha}^{m,n}$ for all $k\geq N_m$. So, $\{x_k\}_{k=N_m}^{\infty}\subset V_{\alpha}^{m,n}$.
Thus, $x\in V_{\alpha}^{m,n}$ as $V_{\alpha}^{m,n}$ is compact. Hence, $x\in\bigcap_{m=0}^{\infty}V_{\alpha}^{m,n}\subset\Lambda_n$.
Therefore, $\Lambda_n$ is closed, and thus compact, since $\bigcup_{i=1}^{N}V_{i,n}$ is compact.

For any $x\in\Lambda_n$, there exists $\alpha\in\Sigma_{N}^{+}(A)$ such that $x\in\bigcap_{m=0}^{\infty}
V_{\alpha}^{m,n}$. Define $\varphi_n(x)=\alpha$. Then, the map $\varphi_n: \Lambda_n\to\Sigma_{N}^{+}(A)$
is well defined since $V_{i,k}\cap V_{j,k}=\emptyset$ for all $1\leq i\neq j\leq N$ and $k\geq0$. Clearly,
$\varphi_n$ is surjective. It is easy to verify that $f_{n}(x)\in\bigcap_{m=0}^{\infty}
V_{\sigma_{A}(\alpha)}^{m,n+1}\subset\Lambda_{n+1}$. Thus, $f_{n}(\Lambda_{n})\subset\Lambda_{n+1}$.
Moreover, $\varphi_{n+1}\circ f_{n}(x)=\sigma_{A}(\alpha)=\sigma_{A}\circ \varphi_{n}(x)$.
Hence, $\varphi_{n+1}\circ f_{n}=\sigma_{A}\circ \varphi_{n}$.

Next, it will be shown that $\varphi_n$ is continuous in $\Lambda_{n}$ for any fixed $n\geq0$. For any $\epsilon>0$,
there exists $N_1\in\mathbf{Z}^{+}$ such that $2^{-N_1}<\epsilon$. Since $\Lambda_{n+j}$ is compact and
$f_{n+j}$ is continuous in $\Lambda_{n+j}$, $0\leq j\leq N_1-1$, there exists $\delta>0$ such that
for any $x,y\in \Lambda_n$ with $\varphi_n(x)=\alpha=(a_0, a_1,\cdots)$ and $\varphi_n(y)=\beta=(b_0, b_1,\cdots)$,
if $d(x,y)<\delta$, then
\vspace{-0.2cm}$$d(f_n^{j}(x),f_n^{j}(y))<\min_{1\leq i_1\neq i_2\leq N}d(V_{i_1,n+j},V_{i_2,n+j}),\;0\leq j\leq N_1.\vspace{-0.2cm}$$
This, together with the fact that $f_{n}^{j}(x)\in V_{a_j,n+j}$
and $f_{n}^{j}(y)\in V_{b_j,n+j}$, $0\leq j\leq N_1$, implies that $a_j=b_j,\;\;0\leq j\leq N_1.$
Thus, one has
\vspace{-0.2cm}$$\hat{\rho}(\varphi_n(x),\varphi_n(y))=\hat{\rho}(\alpha,\beta)\leq2^{-N_1}<\epsilon.\vspace{-0.2cm}$$
Hence, $\varphi_n$ is continuous in $\Lambda_{n}$. Therefore, the invariant subsystem of
system (2.1) on $\{\Lambda_{n}\}_{n=0}^{\infty}$ is topologically semiconjugate
to $(\Sigma_{N}^{+}(A),\sigma_A)$.

{\bf Sufficiency}. Suppose that the invariant subsystem of system (2.1) on
$\{\Lambda_{n}\}_{n=0}^{\infty}$ is topologically $\{\hat{\varphi}_n\}_{n=0}^{\infty}$-semiconjugate
to $(\Sigma_{N}^{+}(A),\sigma_A)$. Then, $\hat{\varphi}_n: \Lambda_n\to\Sigma_{N}^{+}(A)$ is continuous and
surjective, and satisfies that
\vspace{-0.2cm}$$\hat{\varphi}_{n+1}\circ f_n=\sigma_{A}\circ \hat{\varphi}_n,\;n\geq0.\eqno(3.3)\vspace{-0.2cm}$$
Thus,
\vspace{-0.2cm}$$\hat{\varphi}_{n+k}\circ f_n^{k}=\sigma_{A}^{k}\circ\hat{\varphi}_n,\;n\geq0,\;k\geq1.\eqno(3.4)\vspace{-0.2cm}$$
Let
\vspace{-0.2cm}$$V_{i,n}:=\hat{\varphi}_n^{-1}(U_i),\;1\leq i\leq N,\;n\geq0,\eqno(3.5)\vspace{-0.2cm}$$
where $U_i$ is specified in {\rm(2.8)}. For any $1\leq i\leq N$ and $n\geq0$, $V_{i,n}\neq\emptyset$ since $\hat{\varphi}_n$ is surjective;
$V_{i,n}\subset\Lambda_{n}$; and $V_{i,n}$ is compact because $U_i$ is closed, $\varphi_n$ is continuous, and $\Lambda_{n}$ is compact.
It is evident that $E\subset\bigcup_{n=0}^{\infty}\Lambda_{n}$, implying that $f_n$ is continuous in $E$, $n\geq0$.
Since $U_i,\;1\leq i\leq N$, are disjoint, $V_{i,n}\cap V_{j,n}=\emptyset$ for all $1\leq i\neq j\leq N$ and $n\geq0$.
Fix any $m,n\geq0$ and any $\alpha=(a_0,a_1,\cdots)\in\Sigma_{N}^{+}(A)$. It follows from {\rm(3.4)} and {\rm(3.5)} that
\vspace{-0.2cm}$$
V_{\alpha}^{m,n}=\bigcap_{k=0}^{m}f_{n}^{-k}(V_{a_{k},n+k})
=\bigcap_{k=0}^{m}(f_{n}^{-k}\circ\hat{\varphi}_{n+k}^{-1})(U_{a_{k}})
=\bigcap_{k=0}^{m}(\hat{\varphi}_{n}^{-1}\circ\sigma_{A}^{-k})(U_{a_{k}})=\hat{\varphi}_{n}^{-1}\bigg(\bigcap_{k=0}^{m}
\sigma_{A}^{-k}(U_{a_{k}})\bigg).
\vspace{-0.2cm}$$
By Lemma 2.7, $\bigcap_{k=0}^{m}\sigma_{A}^{-k}(U_{a_{k}})\neq\emptyset$.
Hence, $V_{\alpha}^{m,n}\neq\emptyset$, $m,n\geq0$, $\alpha\in\Sigma_{N}^{+}(A)$.
This completes the proof of the theorem.\medskip

Now, a sufficient condition is derived to ensure the system (2.1) to have an invariant subsystem
that is topologically semiconjugate to a subshift of finite type.\medskip

\noindent{\bf Corollary 3.1.} {\it Let assumptions {\rm(i)}-{\rm(ii)} of Lemma 3.1 hold and assume $V_{i,n}\cap V_{j,n}=\emptyset$,
$1\leq i\neq j\leq N$, $n\geq0$. Then, for any $n\geq0$, there exists a nonempty compact subset $\Lambda_{n}\subset \bigcup_{i=1}^{N}V_{i,n}$
with $f_{n}(\Lambda_{n})=\Lambda_{n+1}$ such that $f_n$ is continuous in $\bigcup_{n=0}^{\infty}\Lambda_{n}$ and the invariant subsystem of system {\rm(2.1)}
on $\{\Lambda_{n}\}_{n=0}^{\infty}$ is topologically semiconjugate to $(\Sigma_{N}^{+}(A),\sigma_A)$.}\medskip

\noindent{\bf Proof.} By {\rm(ii)} of Lemma 3.1 and Theorem 3.1, it suffices to show that $\Lambda_{n+1}\subset f_{n}(\Lambda_{n})$
for all $n\geq0$, where $\Lambda_{n}$ is specified in {\rm(3.2)}.  Fix any $n\geq0$. For any $y\in\Lambda_{n+1}$,
there exists $\beta=(b_0,b_1,\cdots)\in\Sigma_{N}^{+}(A)$ such that $y\in\bigcap_{m=0}^{\infty}V_{\beta}^{m,n+1}$
by {\rm(3.2)}. Then, $f_{n+1}^{k}(y)\in V_{b_{k},n+1+k}$ for any $k\geq0$. Since $\sigma_A$ is
surjective, there exists $\beta'=(a_0,b_0,b_1,\cdots)\in\Sigma_{N}^{+}(A)$ such that $\sigma_A(\beta')=\beta$.
As system {\rm(2.1)} is weakly $A$-coupled-expanding in $\{V_{i,n}\}_{n=0}^{\infty}$, $V_{b_{0},n+1}
\subset f_n(V_{a_{0},n})$. Consequently, there exists $x\in V_{a_{0},n}$ such that $y=f_n(x)$. So,
$f_{n}^{k+1}(x)\in V_{b_{k},n+1+k}$, $k\geq0$, and thus $x\in\bigcap_{m=0}^{\infty}V_{\beta'}^{m,n}\subset\Lambda_{n}$.
Hence, $y\in f_{n}(\Lambda_{n})$. Therefore, $\Lambda_{n+1}\subset f_{n}(\Lambda_{n})$ for all $n\geq0$.
The proof is complete.\medskip

Next, a necessary and sufficient condition is established under which an invariant subsystem of system (2.1) is topologically conjugate
to a subshift of finite type.\medskip

\noindent{\bf Theorem 3.2.} {\it Let $A=(a_{ij})_{N\times N}$ be a transition matrix. Then, there exists a sequence
$\{V_{i,n}\}_{n=0}^{\infty}$ of nonempty compact subsets of $X$ with $V_{i,n}\cap V_{j,n}=\emptyset$,
$1\leq i\neq j\leq N$, $n\geq0$, such that $f_n$ is continuous in $E$, $n\geq0$, and $\bigcap_{m=0}^{\infty}
V_{\alpha}^{m,n}$ is a singleton for any $n\geq0$ and any $\alpha\in\Sigma_{N}^{+}(A)$ if and only if, for any $n\geq0$,
there exists a nonempty compact subset $\Lambda_{n}\subset X$ with $f_{n}(\Lambda_{n})=\Lambda_{n+1}$ such that $f_n$
is continuous in $\bigcup_{n=0}^{\infty}\Lambda_{n}$ and the invariant subsystem of system {\rm(2.1)} on
$\{\Lambda_{n}\}_{n=0}^{\infty}$ is topologically conjugate to $(\Sigma_{N}^{+}(A),\sigma_A)$.}\medskip

\noindent{\bf Proof.} {\bf Necessity}. Fix any $n\geq0$. Denote
\vspace{-0.2cm}$$
\bigcap_{m=0}^{\infty}V_{\alpha}^{m,n}:=\{x^{n}(\alpha)\},\;\;\Lambda_{n}:=\bigcup_{\alpha\in\Sigma_{N}^{+}(A)}\{x^{n}(\alpha)\}.
\vspace{-0.2cm}$$
Clearly, $\Lambda_{n}\neq\emptyset$ and $\Lambda_{n}\subset\bigcup_{i=1}^{N}V_{i,n}$. Then, $\bigcup_{n=0}^{\infty}\Lambda_{n}\subset E$,
and thus $f_n$ is continuous in $\bigcup_{n=0}^{\infty}\Lambda_{n}$. It can be easily verified that, for any $\alpha\in\Sigma_{N}^{+}(A)$,
\vspace{-0.2cm}$$
f_{n}(x^{n}(\alpha))=x^{n+1}(\sigma_{A}(\alpha)).
\eqno(3.6)\vspace{-0.2cm}$$
This, together with the fact that $\sigma_{A}$ is surjective, implies that $f_{n}(\Lambda_{n})=\Lambda_{n+1}$.

Define a map $h_{n}:\Sigma_{N}^{+}(A)\to\Lambda_{n}$ by $h_{n}(\alpha)=x^{n}(\alpha)$. Then, $h_n$ is well defined, surjective, and one-to-one,
since $V_{i,k}\cap V_{j,k}=\emptyset$ for all $1\leq i\neq j\leq N$ and $k\geq0$. By {\rm(3.6)}, for any
$\alpha\in\Sigma_{N}^{+}(A)$, one has
\vspace{-0.2cm}$$
h_{n+1}\circ\sigma_{A}(\alpha)=x^{n+1}(\sigma_{A}(\alpha))=f_{n}(x^{n}(\alpha))=f_{n}\circ h_n(\alpha),
\vspace{-0.2cm}$$
which yields that
\vspace{-0.2cm}$$
h_{n+1}\circ\sigma_{A}=f_{n}\circ h_{n}.
\eqno(3.7)\vspace{-0.2cm}$$

Next, it will be shown that $h_n$ is continuous in $\Sigma_{N}^{+}(A)$ for any fixed $n\geq0$. Fix any $\alpha\in\Sigma_{N}^{+}(A)$.
By the assumption that $\bigcap_{m=0}^{\infty}V_{\alpha}^{m,n}$ is a singleton and using Lemma 2.3, one has that $\lim_{m\to\infty}
d(V_{\alpha}^{m,n})=0$. Thus, for any $\epsilon>0$, there exists $N_2\in\mathbf{Z}^{+}$ such that $d(V_{\alpha}^{m,n})<\epsilon$
for all $m\geq N_2$. Denote $\delta:=2^{-N_{2}}$. For any $\beta\in\Sigma_{N}^{+}(A)$ with $\hat{\rho}(\alpha,\beta)<\delta$,
one has that $a_j=b_j,\;\;0\leq j\leq N_2$. Then, $x^{n}(\alpha),x^{n}(\beta)\in V_{\alpha}^{N_2,n}$, and thus
\vspace{-0.2cm}$$d(h_n(\alpha),h_n(\beta))=d(x^{n}(\alpha),x^{n}(\beta))\leq d(V_{\alpha}^{N_2,n})<\epsilon.\vspace{-0.2cm}$$
Therefore, $h_n$ is a continuous and bijective map from a compact metric space $(\Sigma_{N}^{+}(A),\hat{\rho})$
to a metric space $(\Lambda_n,d)$. So, $\Lambda_{n}=h_{n}(\Sigma_{N}^{+}(A))$ is compact and $h_n$ is a homeomorphism.
Consequently, the invariant subsystem of system {\rm(2.1)} on $\{\Lambda_{n}\}_{n=0}^{\infty}$ is topologically conjugate
to $(\Sigma_{N}^{+}(A),\sigma_A)$.

{\bf Sufficiency}. Suppose that the invariant subsystem of system {\rm(2.1)} on
$\{\Lambda_{n}\}_{n=0}^{\infty}$ is topologically $\{\hat{\varphi}_n\}_{n=0}^{\infty}$-conjugate to
$(\Sigma_{N}^{+}(A),\sigma_A)$. Then, $\hat{\varphi}_n: \Lambda_{n}\to\Sigma_{N}^{+}(A)$ is a homeomorphism,
$n\geq0$, and satisfies {\rm(3.3)} and {\rm(3.4)}. Let $V_{i,n}$ be specified in {\rm(3.5)}.
Then, $\{V_{i,n}\}_{n=0}^{\infty}$ is a sequence of nonempty compact subsets of $X$ with $V_{i,n}\cap V_{j,n}=\emptyset$,
$1\leq i\neq j\leq N$, $n\geq0$, and $f_n$ is continuous in $E$, $n\geq0$. Fix any $n\geq0$
and any $\alpha=(a_0,a_1,\cdots)\in\Sigma_{N}^{+}(A)$. Then, by {\rm(3.4)} and {\rm(3.5)}, one has that
\begin{align*}
\bigcap_{m=0}^{\infty}V_{\alpha}^{m,n}=&\bigcap_{m=0}^{\infty}f_{n}^{-m}(V_{a_{m},n+m})
=\bigcap_{m=0}^{\infty}(f_{n}^{-m}\circ\hat{\varphi}_{n+m}^{-1})(U_{a_{m}})\\
=&\bigcap_{m=0}^{\infty}(\hat{\varphi}_{n}^{-1}\circ\sigma_{A}^{-m})(U_{a_{m}})
=\hat{\varphi}_{n}^{-1}\bigg(\bigcap_{m=0}^{\infty}\sigma_{A}^{-m}(U_{a_{m}})\bigg).
\end{align*}
Note that $\bigcap_{m=0}^{\infty}\sigma_{A}^{-m}(U_{a_{m}})=\{\alpha\}$ by Proposition 3.1 in \cite{shicsf}.
Hence, $\bigcap_{m=0}^{\infty}V_{\alpha}^{m,n}$ is a singleton for any $n\geq0$ and any $\alpha\in\Sigma_{N}^{+}(A)$.
This completes the proof of the theorem.\medskip

\noindent{\bf Remark 3.1.} Theorems 3.1 and 3.2 and Corollary 3.1 extend Theorems 4.1, 3.1, and 4.2 in \cite{ju16}
from autonomous discrete systems to non-autonomous discrete systems, respectively.\medskip

The following result is a direct consequence of (iv) of Lemma 3.1 and Theorem 3.2.\medskip

\noindent{\bf Corollary 3.2.} {\it Let assumptions {\rm(i)}-{\rm(iii)} of Lemma 3.1 hold and suppose that
$V_{i,n}\cap V_{j,n}=\emptyset$, $1\leq i\neq j\leq N$, $n\geq0$. Then, for any $n\geq0$,
there exists a nonempty compact subset $\Lambda_{n}\subset\bigcup_{i=1}^{N}V_{i,n}$ with $f_{n}(\Lambda_{n})=\Lambda_{n+1}$ such that
$f_n$ is continuous in $\bigcup_{n=0}^{\infty}\Lambda_{n}$ and the invariant subsystem of system {\rm(2.1)} on $\{\Lambda_{n}\}_{n=0}^{\infty}$
is topologically conjugate to $(\Sigma_{N}^{+}(A),\sigma_A)$.}\medskip

The next result shows that, under certain conditions, topological conjugacy to a subshift of finite type implies weak coupled-expansion
for a transition matrix for system {\rm(2.1)}.\medskip

\noindent{\bf Proposition 3.1.} {\it Let $A=(a_{ij})_{N\times N}$ be a transition matrix. Assume that, for any $n\geq0$, there
exists a nonempty compact subset $\Lambda_{n}\subset X$ with $f_{n}(\Lambda_{n})=\Lambda_{n+1}$ such that $f_n$ is continuous
in $\bigcup_{n=0}^{\infty}\Lambda_{n}$ and the invariant subsystem of system {\rm(2.1)} on $\{\Lambda_{n}\}_{n=0}^{\infty}$
is topologically conjugate to $(\Sigma_{N}^{+}(A),\sigma_A)$. Then, there exists a sequence $\{V_{i,n}\}_{n=0}^{\infty}$
of nonempty compact subsets of $X$, $1\leq i\leq N$, with $\bigcup_{i=1}^{N}V_{i,n}=\Lambda_{n}$, such that
$f_n$ is continuous in $E$, $n\geq0$, and system {\rm(2.1)} is strictly weakly $A$-coupled-expanding in
$\{V_{i,n}\}_{n=0}^{\infty}$, $1\leq i\leq N$.}\medskip

\noindent{\bf Proof.} Suppose that the invariant subsystem of system {\rm(2.1)} on $\{\Lambda_{n}\}_{n=0}^{\infty}$ is
topologically $\{\hat{\varphi}_n\}_{n=0}^{\infty}$-conjugate to $(\Sigma_{N}^{+}(A),\sigma_A)$. Then, $\hat{\varphi}_n: \Lambda_{n}
\to\Sigma_{N}^{+}(A)$ is a homeomorphism, $n\geq0$, and satisfies {\rm(3.3)}. Let $V_{i,n}$ be specified in {\rm(3.5)}.
Then, $f_n$ is continuous in $E$, $n\geq0$, and $\{V_{i,n}\}_{n=0}^{\infty}$ is a sequence of nonempty compact subsets
of $X$ with $V_{i,n}\cap V_{j,n}=\emptyset$, $1\leq i\neq j\leq N$, $n\geq0$, since $U_i,\;1\leq i\leq N$, are disjoint.
Thus, $d(V_{i,n},V_{j,n})>0$, $1\leq i\neq j\leq N$, $n\geq0$. By {\rm(2.9)} and {\rm(3.5)}, one has that
\vspace{-0.2cm}$$\bigcup_{i=1}^{N}V_{i,n}=\bigcup_{i=1}^{N}\hat{\varphi}_{n}^{-1}(U_i)=\hat{\varphi}_{n}^{-1}\bigg(\bigcup_{i=1}^{N}U_i\bigg)
=\hat{\varphi}_{n}^{-1}(\Sigma_{N}^{+}(A))=\Lambda_{n},\;n\geq0.\vspace{-0.2cm}$$
Fix any $n\geq0$ and $1\leq i\leq N$. By {\rm(3.3)}, {\rm(3.5)}, and Lemma 2.7, one has that
\vspace{-0.2cm}$$f_{n}(V_{i,n})=f_{n}(\hat{\varphi}_{n}^{-1}(U_i))=\hat{\varphi}_{n+1}^{-1}(\sigma_{A}(U_i))
\supset\hat{\varphi}_{n+1}^{-1}\bigg(\bigcup_{a_{ij}=1}U_j\bigg)=\bigcup_{a_{ij}=1}\hat{\varphi}_{n+1}^{-1}(U_j)=\bigcup_{a_{ij}=1}V_{j,n+1}.
\vspace{-0.2cm}$$
Consequently, system {\rm(2.1)} is strictly weakly $A$-coupled-expanding in $\{V_{i,n}\}_{n=0}^{\infty}$, $1\leq i\leq N$.
The proof is complete.\bigskip

\noindent{\bf 3.2. Topological equi-semiconjugacy and equi-conjugacy}\medskip

First, it will be shown that, under certain conditions, system (2.1) has an invariant subsystem that is topologically equi-semiconjugate
to a subshift of finite type.\medskip

\noindent{\bf Theorem 3.3.} {\it Let $A=(a_{ij})_{N\times N}$ be a transition matrix and $V_i$, $1\leq i \leq N$,
be nonempty subsets of $X$ with $d(V_i,V_j)>0$, $1\leq i\neq j\leq N$. Assume that there exists a
nonempty compact subset $V_{i,n}$ of $X$ with $V_{i,n}\subset V_i$, $1\leq i\leq N$, $n\geq0$,
such that $\{f_n\}_{n=0}^{\infty}$ is equi-continuous in $E$ and satisfies assumption ${\bf(H_1)}$.
Then, for any $n\geq0$, there exists a nonempty compact subset $\Lambda_{n}\subset\bigcup_{i=1}^{N}V_{i,n}$
with $f_{n}(\Lambda_{n})\subset\Lambda_{n+1}$ such that the invariant subsystem of system {\rm(2.1)}
on $\{\Lambda_{n}\}_{n=0}^{\infty}$ is topologically equi-semiconjugate to $(\Sigma_{N}^{+}(A),\sigma_A)$.
Consequently, $h(f_{0,\infty},\Lambda_0)$ $\geq\log\rho(A)$ in the case that $(X,d)$ is compact,
where $\rho(A)$ is specified in {\rm(2.7)}.}\medskip

\noindent{\bf Proof.} By Theorem 3.1, it suffices to show that $\{\varphi_n\}_{n=0}^{\infty}$ is
equi-continuous in $\{\Lambda_n\}_{n=0}^{\infty}$, where $\Lambda_n$ and $\varphi_n: \Lambda_n\to\Sigma_{N}^{+}(A)$ are
specified in the necessity part of the proof of Theorem 3.1.

For any $\epsilon>0$, there exists $N_3\in\mathbf{Z}^{+}$
such that $2^{-N_3}<\epsilon$. Since $\{f_n\}_{n=0}^{\infty}$ is equi-continuous in $E$,
there exists $\delta>0$ such that, for any $n\geq0$ and any $x,y\in \Lambda_n$ with $\varphi_n(x)=\alpha=
(a_0, a_1,\cdots)$ and $\varphi_n(y)=\beta=(b_0, b_1,\cdots)$, if $d(x,y)<\delta$, then
\vspace{-0.2cm}$$d(f_n^{j}(x),f_n^{j}(y))<\min_{1\leq i\neq j\leq N}d(V_i,V_j),\;0\leq j\leq N_3.\vspace{-0.2cm}$$
This, together with the fact that
$f_{n}^{j}(x)\in V_{a_j}$ and $f_{n}^{j}(y)\in V_{b_j}$, implies that $a_j=b_j$, $0\leq j\leq N_3$.
Thus, one has
\vspace{-0.2cm}$$\hat{\rho}(\varphi_{n}(x),\varphi_{n}(y))=\hat{\rho}(\alpha,\beta)\leq2^{-N_3}<\epsilon.\vspace{-0.2cm}$$
Hence, $\{\varphi_n\}_{n=0}^{\infty}$ is equi-continuous in $\{\Lambda_n\}_{n=0}^{\infty}$.
Therefore, the invariant subsystem of system {\rm(2.1)} on $\{\Lambda_{n}\}_{n=0}^{\infty}$
is topologically equi-semiconjugate to $(\Sigma_{N}^{+}(A),\sigma_A)$, and consequently
$h(f_{0,\infty},\Lambda_0)$ $\geq h(\sigma_A, \Sigma_{N}^{+}(A))=\log\rho(A)$ by Lemma 2.2.
This completes the proof.\medskip

\noindent{\bf Remark 3.2.} The space $X$ is required to be compact in the notion of topological entropy
for non-autonomous discrete systems introduced in \cite{kolyada}. So, this condition is needed
in the present paper.\medskip

Under a stronger and more verifiable condition, the following stronger conclusion can be drawn.\medskip

\noindent{\bf Theorem 3.4.} {\it Let all the assumptions of Theorem 3.3 hold, except that assumption ${\bf(H_1)}$
is replaced by that system {\rm(2.1)} is weakly $A$-coupled-expanding in $\{V_{i,n}\}_{n=0}^{\infty}$, $1\leq i \leq N$.
Then, all the conclusions of Theorem 3.3 hold and $f_{n}(\Lambda_{n})=\Lambda_{n+1}$ for all $n\geq0$.}\medskip

\noindent{\bf Proof.} By (ii) of Lemma 3.1, one can verify that all the assumptions of Theorem 3.3 hold. Thus,
all the conclusions of Theorem 3.3 can be derived. By the same method used in the proof of Corollary 3.1,
one can prove that $f_{n}(\Lambda_{n})=\Lambda_{n+1}$ for all $n\geq0$. The proof is complete.\medskip

Next, some sufficient conditions are derived to ensure a subshift of finite type to be topologically equi-semiconjugate
to an invariant subsystem of system {\rm(2.1)}.\medskip

\noindent{\bf Theorem 3.5.} {\it Let $A=(a_{ij})_{N\times N}$ be a transition matrix. Assume that there exists
a sequence $\{V_{i,n}\}_{n=0}^{\infty}$ of nonempty compact subsets of $X$, $1\leq i\leq N$, such that
$f_n$ is continuous in $E$, $n\geq0$, and satisfies assumptions ${\bf(H_1)}$ and ${\bf(H_2)}$.
Then, for any $n\geq0$, there exists a nonempty compact subset $\Lambda_{n}\subset\bigcup_{i=1}^{N}V_{i,n}$ with
$f_{n}(\Lambda_{n})=\Lambda_{n+1}$ such that $(\Sigma_{N}^{+}(A),\sigma_A)$ is topologically equi-semiconjugate
to the invariant subsystem of system {\rm(2.1)} on $\{\Lambda_{n}\}_{n=0}^{\infty}$. Consequently,
$h(f_{0,\infty},\Lambda_0)\leq\log\rho(A)$ in the case that $(X,d)$ is compact.}\medskip

\noindent{\bf Proof.} By Lemma 2.3 and (i) of Lemma 3.1, one has that $\bigcap_{m=0}^{\infty}V_{\alpha}^{m,n}$
is a singleton, $n\geq0$, $\alpha\in\Sigma_{N}^{+}(A)$. Let $\Lambda_{n}$ and $h_n: \Sigma_{N}^{+}(A)\to\Lambda_{n}$,
$n\geq0$, be specified as in the necessity part of the proof of Theorem 3.2. Then, $h_n$ is well defined and surjective,
and satisfies {\rm(3.7)} for all $n\geq0$.

It suffices to show that $\{h_n\}_{n=0}^{\infty}$ is equi-continuous in $\Sigma_{N}^{+}(A)$.
Fix any $\alpha=(a_0,a_1,\cdots)\in\Sigma_{N}^{+}(A)$. It follows from assumption ${\bf(H_2)}$ that, for any $\epsilon>0$,
there exists $N_4\in\mathbf{Z}^{+}$ such that $d(V_{\alpha}^{m,n})<\epsilon$, $m\geq N_4$, $n\geq0$.
Denote $\delta:=2^{-N_{4}}$. For any $n\geq0$ and any $\beta=(b_0,b_1,\cdots)\in\Sigma_{N}^{+}(A)$ with
$\hat{\rho}(\alpha,\beta)<\delta$, one has that $a_j=b_j$, $0\leq j\leq N_4$. So, $h_n(\alpha),h_n(\beta)
\in V_{\alpha}^{N_4,n}$. Thus
\vspace{-0.2cm}$$d(h_n(\alpha),h_n(\beta))\leq d(V_{\alpha}^{N_4,n})<\epsilon.\vspace{-0.2cm}$$
Hence, $\{h_n\}_{n=0}^{\infty}$ is equi-continuous at $\alpha$. This, together with the fact that
$(\Sigma_{N}^{+}(A),\hat{\rho})$ is compact, implies that $\{h_n\}_{n=0}^{\infty}$ is equi-continuous
in $\Sigma_{N}^{+}(A)$ by Lemma 2.5 in \cite{Shao16}. Therefore, $(\Sigma_{N}^{+}(A),\sigma_A)$ is topologically
equi-semiconjugate to the invariant subsystem of system {\rm(2.1)} on $\{\Lambda_{n}\}_{n=0}^{\infty}$.
Consequently, $h(f_{0,\infty},\Lambda_0)\leq h(\sigma_A,\Sigma_{N}^{+}(A))=\log\rho(A)$ by Lemma 2.2.
This completes the proof.\medskip

The following result can be directly derived from Lemma 3.1 and Theorem 3.5.\medskip

\noindent{\bf Corollary 3.3.} {\it Let assumptions {\rm(i)}-{\rm(iii)} of Lemma 3.1 hold.
Then, all the conclusions of Theorem 3.5 hold.}\medskip

Based on Theorems 3.3 and 3.5, it can be shown that, under certain conditions, system (2.1) has an invariant subsystem that is
topologically equi-conjugate to a subshift of finite type.\medskip

\noindent{\bf Theorem 3.6.} {\it Let all the assumptions of Theorem 3.3 and assumption ${\bf(H_2)}$ hold.
Then, for any $n\geq0$, there exists a nonempty compact subset $\Lambda_{n}\subset\bigcup_{i=1}^{N}V_{i,n}$
with $f_{n}(\Lambda_{n})=\Lambda_{n+1}$ such that the invariant subsystem of system {\rm(2.1)} on
$\{\Lambda_{n}\}_{n=0}^{\infty}$ is topologically equi-conjugate to $(\Sigma_{N}^{+}(A),\sigma_A)$.
Consequently, $h(f_{0,\infty},\Lambda_0)=\log\rho(A)$ in the case that $(X,d)$ is compact.
Further, if $A$ is irreducible and $\sum_{j=1}^{N}a_{i_{0}j}\geq2$ for some $1\leq i_0\leq N$,
then system {\rm(2.1)} is Li-Yorke $\delta$-chaotic for some $\delta>0$, and is also
distributionally chaotic.}\medskip

\noindent{\bf Proof.} Let $\varphi_{n}: \Lambda_{n}\to\Sigma_{N}^{+}(A)$ and $h_{n}:\Sigma_{N}^{+}(A)\to\Lambda_{n}$,
$n\geq0$, be specified as in the proofs of Theorems 3.3 and 3.5, respectively. Then, $\varphi_n=h_n^{-1}$, $n\geq0$,
are homeomorphisms, $\{\varphi_{n}\}_{n=0}^{\infty}$ and $\{h_{n}\}_{n=0}^{\infty}$ are equi-continuous in
$\{\Lambda_{n}\}_{n=0}^{\infty}$ and $\Sigma_{N}^{+}(A)$, respectively, and {\rm(3.7)} holds
for all $n\geq0$. Hence, the invariant subsystem of system {\rm(2.1)} on $\{\Lambda_{n}\}_{n=0}^{\infty}$ is topologically
equi-conjugate to $(\Sigma_{N}^{+}(A),\sigma_A)$, and consequently $h(f_{0,\infty},\Lambda_0)=\log\rho(A)$ by Lemma 2.2.

If $A$ is irreducible and $\sum_{j=1}^{N}a_{i_{0}j}\geq2$ for some $1\leq i_0\leq N$, then $(\Sigma_{N}^{+}(A),\sigma_{A})$
has an uncountable $\gamma$-scrambled set with $\gamma=1$ by Lemma 4.1 in \cite{Shi09}. This, together with the result of
(i) in Lemma 2.1, implies that system {\rm(2.1)} is Li-Yorke $\delta$-chaotic for some $\delta>0$. It is also distributionally
chaotic by (ii) in Lemma 2.1, since $(\Sigma_{N}^{+}(A),\sigma_{A})$ is distributionally chaotic by Theorem 1.4 in \cite{wanghy}.
The proof is complete.\medskip

%The following result is a direct consequence of (ii)-(iii) of Lemma 3.1 and Theorem 3.6.\medskip

%By Lemma 3.1 and Theorem 3.6 one can get the following stronger conclusion.\medskip

One can obtain the following stronger result under some more verifiable conditions.\medskip

\noindent{\bf Theorem 3.7.} {\it Let all the assumptions of Theorem 3.4 and assumption {$\rm(iii)$} of Lemma 3.1 hold.
Then, all the conclusions of Theorem 3.6 hold and system {\rm(2.1)} is chaotic in the strong sense of Li-Yorke
in the case that $A$ is irreducible and $\sum_{j=1}^{N}a_{i_{0}j}\geq2$ for some $1\leq i_0\leq N$.}\medskip

\noindent{\bf Proof.} It follows from (ii)-(iii) in Lemma 3.1 that all the assumptions of Theorem 3.6 hold, thus
all the conclusions of Theorem 3.6 hold. By assumptions of {$\rm(iii)$} in Lemma 3.1 that $V_i$, $1\leq i\leq N$,
are bounded, one has that $\cup_{n=0}^{\infty}\Lambda_{n}$ is bounded, since $\cup_{n=0}^{\infty}\Lambda_{n}\subset\bigcup_{i=1}^{N}V_{i,n}\subset\bigcup_{i=1}^{N}V_{i}$. This implies that
system {\rm(2.1)} is chaotic in the strong sense of Li-Yorke.\medskip

By the same method as that used in the proof of Theorem 3.1 in \cite{Shao16}, one can get the following result.\medskip

\noindent{\bf Theorem 3.8.} {\it Let assumptions {\rm(i)}-{\rm(ii)} of Lemma 3.1 hold, $f_n$ be continuous in the
compact metric space $(X,d)$, $n\geq0$, and $V_i$, $1\leq i\leq N$, be disjoint nonempty closed subsets of $X$ with
$V_{i,n}\subset V_i$ for all $n\geq0$. Then, $h(f_{0,\infty})\geq\log\rho(A)$.}\medskip

\noindent{\bf Remark 3.3.}
\begin{itemize}\vspace{-0.2cm}
\item[{\rm (i)}] Theorems 3.4, 3.5, 3.7, and 3.8 and Corollary 3.3 generalize Theorems 4.4, 4.1, 4.5, 3.1, and 4.2 in \cite{Shao16},
respectively, where only a special case that $V_{i,n}=V_i$, $1\leq i\leq N$, $n\geq0$, is considered.\vspace{-0.2cm}
\item[{\rm (ii)}] Theorems 3.4 and 3.7 relax the assumptions of Theorems 4.4 and 4.5 in \cite{Shao16}, respectively, since
assumption $({\rm ii})_b$ of Theorem 4.4 (resp. Theorem 4.5) in \cite{Shao16} is replaced by a weaker one that
$\{f_n\}_{n=0}^{\infty}$ is equi-continuous in $E$ in Theorems 3.4 and 3.7 here.\vspace{-0.2cm}
\item[{\rm (iii)}] Theorem 3.5 and Corollary 3.3 also relax the assumptions of Theorems 4.1 and 4.2 in \cite{Shao16}, respectively,
since it is only required $f_n$ be continuous in $E$, $n\geq0$, in Theorem 3.5 and Corollary 3.3 here.\vspace{-0.2cm}
\item[{\rm (iv)}] Assumption ${\bf(H_1)}$ in Theorem 3.5 is strictly weaker than assumption (i) of Theorem 4.1 in
\cite{Shao16}, since the converse of {\rm(ii)} of Lemma 3.1 is not true in general, even for autonomous
dynamical systems (see Example 3.1.1 in \cite{ju17}).\vspace{-0.2cm}
\end{itemize}\vspace{-0.2cm}

\bigskip

\noindent{\bf 4. Topological (semi)conjugacy and equi-(semi)conjugacy between induced set-valued systems
and subshifts of finite type}\medskip

First, some relationships of several related dynamical behaviors between system {\rm(2.1)} and its induced set-valued system
{\rm(2.6)} are established in Section 4.1. Then, by these results, together with the results obtained in Section 3,
the topological (semi)conjugacy and equi-(semi)conjugacy between system {\rm(2.6)} and a subshift of finite type
are proved in Sections 4.2 and 4.3, respectively.\medskip

\noindent{\bf 4.1. Relationships of some related dynamical behaviors}\medskip

To start, the following assumption is made.\medskip

${\bf(H_3)}$ Let $(X,d)$ be a compact metric space, $f_n$ be continuous in $X$, $n\geq0$,
$\{V_{i,n}\}_{n=0}^{\infty}$ be a sequence of nonempty closed subsets of $X$, $1\leq i\leq N$,
and $A$ be a transition matrix.\medskip

The following result can be easily verified based on Lemma 2.4 and (i) of Lemma 2.5.\medskip

\noindent{\bf Proposition 4.1.} {\it Let assumption ${\bf(H_3)}$ hold. Then, $\bar{f}_n$ is continuous in the
compact metric space $(\mathcal{K}(X),\mathcal{H}_{d})$, $n\geq0$, and $\{\langle V_{i,n}\rangle\}_{n=0}^{\infty}$
is a sequence of nonempty compact subsets of $\mathcal{K}(X)$, $1\leq i\leq N$. In addition, $V_{i,n}\cap V_{j,n}=\emptyset$ if and only if
$\langle V_{i,n}\rangle\cap\langle V_{j,n}\rangle=\emptyset$
for all $1\leq i\neq j\leq N$ and $n\geq0$.}\medskip

\noindent{\bf Proposition 4.2.} {\it Let assumption ${\bf(H_3)}$ hold. Fix any $m,n\geq0$ and any $\alpha=(a_0,a_1,\cdots)$
$\in\Sigma_{N}^{+}(A)$. Then, $V_{\alpha}^{m,n}\neq\emptyset$ if and only if $\mathcal{H}_{\alpha}^{m,n}\neq\emptyset$, where
\vspace{-0.2cm}$$\mathcal{H}_{\alpha}^{m,n}:=\bigcap_{k=0}^{m}\bar{f}_{n}^{-k}(\langle V_{a_{k},n+k}\rangle).               \eqno(4.1)\vspace{-0.2cm}$$}\medskip

\noindent{\bf Proof.} Suppose that $V_{\alpha}^{m,n}\neq\emptyset$. It follows from Lemma 2.6 and (ii) of Lemma 2.5 that
\vspace{-0.2cm}$$
\mathcal{H}_{\alpha}^{m,n}=\bigcap_{k=0}^{m}\bar{f}_{n}^{-k}(\langle V_{a_{k},n+k}\rangle)\supset\bigcap_{k=0}^{m}
\langle f_{n}^{-k}(V_{a_{k},n+k})\rangle=\bigg\langle\bigcap_{k=0}^{m}f_{n}^{-k}(V_{a_{k},n+k})\bigg\rangle=\langle V_{\alpha}^{m,n}\rangle
\neq\emptyset.
\eqno(4.2)\vspace{-0.2cm}$$

Conversely, suppose that $\mathcal{H}_{\alpha}^{m,n}\neq\emptyset$. Then, there exists
$K_0\in\mathcal{H}_{\alpha}^{m,n}$. So, $\bar{f}_{n}^{k}(K_0)=f_{n}^{k}(K_0)\in\langle V_{a_{k},n+k}\rangle$, and thus
$f_{n}^{k}(K_0)\subset V_{a_{k},n+k}$, $0\leq k\leq m$. Hence,
\vspace{-0.2cm}$$K_0\subset\bigcap_{k=0}^{m}f_{n}^{-k}(V_{a_{k},n+k})=V_{\alpha}^{m,n}.\vspace{-0.2cm}$$
Therefore, $V_{\alpha}^{m,n}\neq\emptyset$. The proof is complete.\medskip

\noindent{\bf Proposition 4.3.} {\it Let assumption ${\bf(H_3)}$ hold. Then, $d(V_{\alpha}^{m,n})$ converges to $0$ as $m\to\infty$
if and only if $\mathcal{H}_{d}(\mathcal{H}_{\alpha}^{m,n})$ converges to $0$ as $m\to\infty$ for any $\alpha=(a_0,a_1,\cdots)
\in\Sigma_{N}^{+}(A)$ and $n\geq0$. Further, $d(V_{\alpha}^{m,n})$ uniformly converges to $0$ with respect to $n\geq0$ as $m\to\infty$
if and only if $\mathcal{H}_{d}(\mathcal{H}_{\alpha}^{m,n})$ uniformly converges to $0$ with respect to $n\geq0$ as $m\to\infty$, for any
$\alpha\in\Sigma_{N}^{+}(A)$.}\medskip

\noindent{\bf Proof.} Fix any $\alpha\in\Sigma_{N}^{+}(A)$ and $n\geq0$. Suppose that $d(V_{\alpha}^{m,n})$ converges to $0$ as $m\to\infty$.
Then, for any $\epsilon>0$, there exists $M_1>0$ such that $d(V_{\alpha}^{m,n})<\epsilon$ for all $m\geq M_1$. Fix any $m\geq M_1$.
For any $K_1, K_2\in \mathcal{H}_{\alpha}^{m,n}$, one has that
\vspace{-0.2cm}$$f_{n}^{k}(K_1),f_{n}^{k}(K_2)\subset V_{a_{k},n+k},\;0\leq k\leq m.                                                 \eqno(4.3)\vspace{-0.2cm}$$
For any $a\in K_1$, there exists $b_1\in K_2$ such that $d(a,K_2)=d(a,b_1)$. By (4.3), one has that
$f_{n}^{k}(a),f_{n}^{k}(b_1)\in V_{a_{k},n+k}$ for all $0\leq k\leq m$, hence $a,b_1\in V_{\alpha}^{m,n}$. Thus,
\vspace{-0.2cm}$$d(a,K_2)=d(a,b_1)\leq d(V_{\alpha}^{m,n})<\epsilon.\vspace{-0.2cm}$$
So, $\sup_{a\in K_1}d(a,K_2)\leq\epsilon$. Similarly, one can verify that $\sup_{b\in K_2}d(b,K_1)\leq\epsilon$. This,
together with (2.4), implies that $\mathcal{H}_{d}(K_1,K_2)\leq\epsilon$. Hence, $\mathcal{H}_{d}(\mathcal{H}_{\alpha}^{m,n})
\leq\epsilon$ for all $m\geq M_1$. Consequently, $\mathcal{H}_{d}(\mathcal{H}_{\alpha}^{m,n})$ converges to $0$ as $m\to\infty$.

Conversely, suppose that $\mathcal{H}_{d}(\mathcal{H}_{\alpha}^{m,n})$ converges to $0$ as $m\to\infty$. Then, for any $\epsilon>0$,
there exists $M_2>0$ such that $\mathcal{H}_{d}(\mathcal{H}_{\alpha}^{m,n})<\epsilon$ for all $m\geq M_2$. Fix any $m\geq M_2$.
For any $x,y\in V_{\alpha}^{m,n}$, one has that $f_{n}^{k}(x),f_{n}^{k}(y)\in V_{a_{k},n+k}$, $0\leq k\leq m$. Thus, $\bar{f}_{n}^{k}(\{x\}),\bar{f}_{n}^{k}(\{y\})\in\langle V_{a_{k},n+k}\rangle$, $0\leq k\leq m$,
which implies that $\{x\},\{y\}\in\bigcap_{k=0}^{m}\bar{f}_{n}^{-k}(\langle V_{a_{k},n+k}\rangle)=\mathcal{H}_{\alpha}^{m,n}$. So,
\vspace{-0.2cm}$$d(x,y)=\mathcal{H}_{d}(\{x\},\{y\})\leq\mathcal{H}_{d}(\mathcal{H}_{\alpha}^{m,n})<\epsilon.\vspace{-0.2cm}$$
Hence, $d(V_{\alpha}^{m,n})\leq\epsilon$ for all $m\geq M_2$. Consequently, $d(V_{\alpha}^{m,n})$ converges to $0$ as $m\to\infty$.

Similarly, one can show that $d(V_{\alpha}^{m,n})$ uniformly converges to $0$ with respect to $n\geq0$ as $m\to\infty$ if and only if $\mathcal{H}_{d}(\mathcal{H}_{\alpha}^{m,n})$ uniformly converges to $0$ with respect to $n\geq0$ as $m\to\infty$, for any
$\alpha\in\Sigma_{N}^{+}(A)$. The proof is complete.\medskip

Next, it will be shown that the (strictly) weak $A$-coupled-expansion of system {\rm(2.1)} is equivalent to that of the
induced set-valued system (2.6).\medskip

\noindent{\bf Proposition 4.4.} {\it Let assumption ${\bf(H_3)}$ hold. System {\rm(2.1)} is {\rm(}strictly{\rm)} weakly $A$-coupled-expanding in
$\{ V_{i,n}\}_{n=0}^{\infty}$, $1\leq i\leq N$, if and only if system {\rm(2.6)} is {\rm(}strictly{\rm)}
weakly $A$-coupled-expanding in $\{\langle V_{i,n}\rangle\}_{n=0}^{\infty}$, $1\leq i\leq N$.}\medskip

\noindent{\bf Proof.} Suppose that system {\rm(2.1)} is weakly $A$-coupled-expanding in $\{ V_{i,n}\}_{n=0}^{\infty}$, $1\leq i\leq N$.
Fix any $1\leq i\leq N$ and $n\geq0$. For any $1\leq j\leq N$ with $a_{ij}=1$ and any $K\in\langle V_{j,n+1}\rangle$, one has that
$K\in\mathcal{K}(X)$ and $K\subset V_{j,n+1}\subset f_{n}(V_{i,n})$. Let $K_0:=f_n^{-1}(K)\cap V_{i,n}$.
Then, $K_0\in\mathcal{K}(X)$, since $f_n$ is continuous and both $K$ and $V_{i,n}$ are compact. Thus, $K_0\in\langle V_{i,n}\rangle$.
It is easy to verify that $f_n(K_0)=K$. Thus, $K\in\bar{f}_{n}(\langle V_{i,n}\rangle)$.
Hence, $\langle V_{j,n+1}\rangle\subset\bar{f}_{n}(\langle V_{i,n}\rangle)$.
Therefore, system {\rm(2.6)} is weakly $A$-coupled-expanding in $\{\langle V_{i,n}\rangle\}_{n=0}^{\infty}$, $1\leq i\leq N$.

Conversely, suppose that  system {\rm(2.6)} is  weakly $A$-coupled-expanding in $\{\langle V_{i,n}\rangle\}_{n=0}^{\infty}$,
$1\leq i\leq N$. Fix any $1\leq i\leq N$ and $n\geq0$. For any $1\leq j\leq N$ with $a_{ij}=1$ and any $y\in V_{j,n+1}$, one has that
$\{y\}\in \langle V_{j,n+1}\rangle\subset\bar{f}_{n}(\langle V_{i,n}\rangle)$. Then, there exists $K_1\in \langle V_{i,n}\rangle$
such that $\{y\}=\bar{f}_{n}(K_1)=f_n(K_1)$. Thus, there exists $x\in K_1\subset V_{i,n}$ such that $y=f_n(x)\in f_n(V_{i,n})$.
Hence, $V_{j,n+1}\subset f_{n}(V_{i,n})$.  Therefore, system {\rm(2.1)} is weakly $A$-coupled-expanding in
$\{V_{i,n}\}_{n=0}^{\infty}$, $1\leq i\leq N$.

It is easy to verify that $d(V_{i,n}, V_{j,n})>0$ if and only if $\mathcal{H}_{d}(\langle V_{i,n}\rangle,\langle V_{j,n}\rangle)>0$,
$1\leq i\neq j\leq N$, $n\geq0$, by Proposition 4.1. Hence, system {\rm(2.1)} is strictly weakly $A$-coupled-expanding
in $\{V_{i,n}\}_{n=0}^{\infty}$, $1\leq i\leq N$, if and only if system {\rm(2.6)} is  strictly weakly $A$-coupled-expanding in
$\{\langle V_{i,n}\rangle\}_{n=0}^{\infty}$, $1\leq i\leq N$. This completes the proof.\bigskip

\noindent{\bf 4.2. Topological semi-conjugacy and conjugacy}\medskip

To start, the following assumption is made.\medskip

${\bf(H_4)}$ Let assumption ${\bf(H_3)}$ hold and $V_{i,n}\cap V_{j,n}=\emptyset$ for all $1\leq i\neq j\leq N$ and $n\geq0$.\medskip

First, a sufficient condition is derived for an invariant subsystem of system (2.6) to be topologically semiconjugate
to $(\Sigma_{N}^{+}(A),\sigma_A)$.\medskip

\noindent{\bf Theorem 4.1.} {\it Let assumptions ${\bf(H_1)}$ and ${\bf(H_4)}$ hold. Then, for any $n\geq0$,
there exists a nonempty compact subset $\mathcal{C}_{n}\subset\bigcup_{i=1}^{N}\langle V_{i,n}\rangle$
with $\bar{f}_{n}(\mathcal{C}_{n})\subset\mathcal{C}_{n+1}$ such that the invariant subsystem of system
{\rm(2.6)} on $\{\mathcal{C}_{n}\}_{n=0}^{\infty}$ is topologically semiconjugate to $(\Sigma_{N}^{+}(A),\sigma_A)$.}\medskip

\noindent{\bf Proof.}
It follows from Proposition 4.1 that $\bar{f}_n$ is continuous in $\mathcal{K}(X)$, $n\geq0$, and
$\{\langle V_{i,n}\rangle\}_{n=0}^{\infty}$ is a sequence of nonempty compact subsets of $\mathcal{K}(X)$,
$1\leq i\leq N$, with $\langle V_{i,n}\rangle\cap \langle V_{j,n}\rangle=\emptyset$ for all
$1\leq i\neq j\leq N$ and  $n\geq0$. By Proposition 4.2, one has that $\mathcal{H}_{\alpha}^{m,n}
\neq\emptyset$ for all $m,n\geq0$ and any $\alpha\in\Sigma_{N}^{+}(A)$. Hence, all the assumptions
of Theorem 3.1 hold for system (2.6). Therefore, the conclusion holds by Theorem 3.1.
The proof is complete.\medskip

\noindent{\bf Remark 4.1.}
\begin{itemize}\vspace{-0.2cm}
\item[{\rm (i)}] By {\rm(3.2)}, {\rm(4.2)}, and (ii) of Lemma 2.5, one has that
 \vspace{-0.2cm}$$\mathcal{C}_{n}=\bigcup_{\alpha\in\Sigma_{N}^{+}(A)}\bigcap_{m=0}^{\infty}\mathcal{H}_{\alpha}^{m,n}
\supset\bigcup_{\alpha\in\Sigma_{N}^{+}(A)}\bigcap_{m=0}^{\infty}\langle V_{\alpha}^{m,n}\rangle
=\bigcup_{\alpha\in\Sigma_{N}^{+}(A)}\bigg\langle\bigcap_{m=0}^{\infty}V_{\alpha}^{m,n}\bigg\rangle.\vspace{-0.2cm}$$
\item[{\rm (ii)}] In the special case that $f_n=f$ and $V_{i,n}=V_i$, $1\leq i\leq N$, $n\geq0$, Theorem 4.1 is the
same as Theorem 3.1.5 in \cite{ju17} for autonomous discrete systems.
\end{itemize}\vspace{-0.2cm}

One can obtain the following stronger conclusion under some stronger and more verifiable condition.\medskip

\noindent{\bf Theorem 4.2.} {\it Let assumption ${\bf(H_4)}$ hold. Assume that system {\rm(2.1)} is
weakly $A$-coupled-expanding in $\{ V_{i,n}\}_{n=0}^{\infty}$, $1\leq i\leq N$. Then, the conclusion of Theorem 4.1
holds and $\bar{f}_{n}(\mathcal{C}_{n})=\mathcal{C}_{n+1}$  for all $n\geq0$.}\medskip

\noindent{\bf Proof.}
It follows from Proposition 4.4 that system {\rm(2.6)} is weakly $A$-coupled-expanding in
$\{\langle V_{i,n}\rangle\}_{n=0}^{\infty}$, $1\leq i\leq N$. This, together with the result of Proposition 4.1,
implies that all the assumptions of Corollary 3.1 hold for system {\rm(2.6)}. Hence, the conclusions hold
by Corollary 3.1. This completes the proof.\medskip

Next, it will be shown that, under certain conditions, system {\rm(2.6)} has an invariant subsystem that is topologically conjugate to
$(\Sigma_{N}^{+}(A),\sigma_A)$.\medskip

\noindent{\bf Theorem 4.3.} {\it Let assumptions ${\bf(H_1)}$ and ${\bf(H_4)}$ hold. Assume that $d(V_{\alpha}^{m,n})$
converges to $0$ as $m\to\infty$ for any $n\geq0$ and any $\alpha\in\Sigma_{N}^{+}(A)$. Then, for any $n\geq0$, there
exists a nonempty compact subset $\mathcal{C}_{n}\subset\bigcup_{i=1}^{N}\langle V_{i,n}\rangle$ with $\bar{f}_{n}
(\mathcal{C}_{n})=\mathcal{C}_{n+1}$ such that the invariant subsystem of system {\rm(2.6)} on
$\{\mathcal{C}_{n}\}_{n=0}^{\infty}$ is topologically conjugate to $(\Sigma_{N}^{+}(A),\sigma_A)$.}\medskip

\noindent{\bf Proof.} Fix any $\alpha\in\Sigma_{N}^{+}(A)$ and $n\geq0$. By (4.1) and Propositions 4.1-4.2,
one has that $\mathcal{H}_{\alpha}^{m,n}$ is a nonempty compact subset of $\mathcal{K}(X)$ and satisfies that
$\mathcal{H}_{\alpha}^{m+1,n}\subset\mathcal{H}_{\alpha}^{m,n}$ for all $m\geq0$. It follows from Proposition 4.3
that $\mathcal{H}_{d}(\mathcal{H}_{\alpha}^{m,n})$ converges to $0$ as $m\to\infty$. Applying Lemma 2.3, one obtains
that $\bigcap_{m=0}^{\infty}\mathcal{H}_{\alpha}^{m,n}$ is a singleton for any fixed $\alpha\in\Sigma_{N}^{+}(A)$ and $n\geq0$.
This, together with the result of Proposition 4.1, implies that all the assumptions of Theorem 3.2 hold for system {\rm(2.6)}.
Therefore, the conclusion holds by Theorem 3.2. The proof is complete.\medskip

By Lemma 3.1 and Theorem 4.3, one has the following result.\medskip

\noindent{\bf Corollary 4.1.} {\it Let assumptions {\rm(ii)}-{\rm(iii)} of Lemma 3.1 and assumption ${\bf(H_4)}$ hold.
Then, the conclusion of Theorem 4.3 holds.}\medskip

\noindent{\bf 4.3. Topological equi-semiconjugacy and equi-conjugacy}\medskip

First, some sufficient conditions are established to ensure the system (2.6) to have an invariant subsystem
that is topologically equi-semiconjugate to a subshift of finite type.\medskip

\noindent{\bf Theorem 4.4.} {\it Let assumptions ${\bf(H_1)}$ and ${\bf(H_3)}$ hold, $\{f_n\}_{n=0}^{\infty}$ be equi-continuous
in $X$, and $V_i$, $1\leq i\leq N$, be disjoint closed subsets of $X$ with $V_{i,n}\subset V_i$, $1\leq i\leq N$, $n\geq0$.
Then, for any $n\geq0$, there exists a nonempty compact subset $\mathcal{C}_{n}\subset\bigcup_{i=1}^{N}\langle V_{i,n}\rangle$
with $\bar{f}_{n}(\mathcal{C}_{n})\subset\mathcal{C}_{n+1}$ such that the invariant subsystem of system {\rm(2.6)}
on $\{\mathcal{C}_{n}\}_{n=0}^{\infty}$ is topologically equi-semiconjugate to $(\Sigma_{N}^{+}(A),\sigma_A)$.
Consequently, $h(\bar{f}_{0,\infty},\mathcal{C}_{0})\geq\log\rho(A)$.}\medskip

\noindent{\bf Proof.}
By (i) of Lemma 2.5, one has that $\langle V_{i}\rangle$ and $\langle V_{i,n}\rangle$, $1\leq i\leq N$, $n\geq0$,
are nonempty compact subsets of $\mathcal{K}(X)$ with $\langle V_{i,n}\rangle\subset\langle V_{i}\rangle$.
Since $V_i$, $1\leq i\leq N$, are disjoint, $\langle V_{i}\rangle$, $1\leq i\leq N$, are disjoint.
Thus, $\mathcal{H}_{d}(\langle V_{i}\rangle,\langle V_{j}\rangle)>0$ for all $1\leq i\neq j\leq N$. It follows
from Lemma 2.4 that $\{\bar{f}_n\}_{n=0}^{\infty}$ is equi-continuous in $\mathcal{K}(X)$. In addition,
$\mathcal{H}_{\alpha}^{m,n}\neq\emptyset$ for all $m,n\geq0$ and $\alpha\in\Sigma_{N}^{+}(A)$
by Proposition 4.2. Hence, all the assumptions of Theorem 3.3 hold for system {\rm(2.6)}.
Therefore, all the conclusions hold by Theorem 3.3.
This completes the proof.\medskip

\noindent{\bf Theorem 4.5.} {\it Let all the assumptions of Theorem 4.4 hold, except that assumption ${\bf(H_1)}$
is replaced by that system {\rm(2.1)} is weakly $A$-coupled-expanding in $\{ V_{i,n}\}_{n=0}^{\infty}$, $1\leq i\leq N$.
Then, all the conclusions of Theorem 4.4 hold and $\bar{f}_{n}(\mathcal{C}_{n})=\mathcal{C}_{n+1}$ for all $n\geq0$.}\medskip

\noindent{\bf Proof.} By the method used in the proof of Theorem 4.4, together with the result of Proposition 4.4,
one can show that all the assumptions of Theorem 3.4 hold for system {\rm(2.6)}. Therefore, the conclusions hold by Theorem 3.4.
The proof is complete.\medskip

Now, a sufficient condition is derived under which $(\Sigma_{N}^{+}(A),\sigma_A)$ is topologically equi-semiconjugate
to an invariant subsystem of system {\rm(2.6)}.\medskip

\noindent{\bf Theorem 4.6.} {\it Let assumptions ${\bf(H_1)}$-${\bf(H_3)}$ hold. Then, for any $n\geq0$, there exists
a nonempty compact subset $\mathcal{C}_{n}\subset\bigcup_{i=1}^{N}\langle V_{i,n}\rangle$ with $\bar{f}_{n}(\mathcal{C}_{n})
=\mathcal{C}_{n+1}$ such that $(\Sigma_{N}^{+}(A),\sigma_A)$ is topologically equi-semiconjugate to the invariant subsystem
of system {\rm(2.6)} on $\{\mathcal{C}_{n}\}_{n=0}^{\infty}$. Consequently, $h(\bar{f}_{0,\infty},\mathcal{C}_{0})\leq\log\rho(A)$.}\medskip

\noindent{\bf Proof.} By Propositions 4.1-4.3, all the assumptions of Theorem 3.5 hold for system {\rm(2.6)}.
Therefore, the conclusions hold by Theorem 3.5. This completes the proof.\medskip

The following result is a direct consequence of Lemma 3.1 and Theorem 4.6.\medskip

\noindent{\bf Corollary 4.2.} {\it Let assumptions ${\rm(ii)}$-${\rm(iii)}$ of Lemma 3.1 and assumption ${\bf(H_3)}$ hold.
Then, all the conclusions of Theorem 4.6 hold.}\medskip

Next, some sufficient conditions are derived to ensure the system (2.6) to have an invariant subsystem
that is topologically equi-conjugate to $(\Sigma_{N}^{+}(A),\sigma_A)$.\medskip

\noindent{\bf Theorem 4.7.} {\it Let all the assumptions of Theorem 4.4 and assumption ${\bf(H_2)}$ hold.
Then, for any $n\geq0$, there exists a nonempty compact subset $\mathcal{C}_{n}\subset\bigcup_{i=1}^{N}
\langle V_{i,n}\rangle$ with $\bar{f}_{n}(\mathcal{C}_{n})=\mathcal{C}_{n+1}$ such that the invariant
subsystem of system {\rm(2.6)} on $\{\mathcal{C}_{n}\}_{n=0}^{\infty}$ is topologically equi-conjugate to
$(\Sigma_{N}^{+}(A),\sigma_A)$, and consequently $h(\bar{f}_{0,\infty},\mathcal{C}_0)=\log\rho(A)$.
Further, if $A$ is irreducible and $\sum_{j=1}^{N}a_{i_{0}j}\geq2$ for some $1\leq i_0\leq N$,
then system {\rm(2.6)} is chaotic in the strong sense of Li-Yorke, and is also distributionally
chaotic.}\medskip

\noindent{\bf Proof.}
By the method used in the proof of of Theorem 4.4, together with the result of Proposition 4.3,
one can verify that all the assumptions of Theorem 3.6 hold for system {\rm(2.6)}. Therefore,
all the conclusions of Theorem 3.6 hold. Further, $\cup_{n=0}^{\infty}\mathcal{C}_{n}$ is
bounded since $\cup_{n=0}^{\infty}\mathcal{C}_{n}\subset\cup_{i=1}^{N}\langle V_{i}\rangle$.
This implies that system {\rm(2.6)} is chaotic in the strong sense of Li-Yorke,
completing the proof.\medskip

The following result is a direct consequence of Lemma 3.1 and Theorem 4.7.\medskip

\noindent{\bf Theorem 4.8.} {\it Let all the assumptions of Theorem 4.5 and assumption {\rm(iii)} of Lemma 3.1 hold.
Then all the conclusions of Theorem 4.7 hold and system {\rm(2.6)} is chaotic in the strong sense of Li-Yorke
in the case that $A$ is irreducible and $\sum_{j=1}^{N}a_{i_{0}j}\geq2$ for some $1\leq i_0\leq N$.}\medskip

\noindent{\bf Proof.} It follows from (ii)-(iii) in Lemma 3.1 that all the assumptions of Theorem 4.7 hold, thus
all the conclusions of Theorem 4.7 hold. By the assumption that $V_i$, $1\leq i\leq N$, are bounded, one has that $\Lambda_{n}$
is bounded, since $\Lambda_{n}\subset\bigcup_{i=1}^{N}V_{i,n}\subset\bigcup_{i=1}^{N}V_{i}$. This implies that
system {\rm(2.1)} is chaotic in the strong sense of Li-Yorke.\medskip

\noindent{\bf Theorem 4.9.} {\it Let all the assumptions of Theorem 3.8 hold. Then, $h(\bar{f}_{0,\infty})\geq\log\rho(A)$.}\medskip

\noindent{\bf Proof.} By (i) of Lemma 2.5, one has that $\langle V_{i}\rangle$, $1\leq i\leq N$, are disjoint
nonempty compact subsets of $\mathcal{K}(X)$ with $\langle V_{i,n}\rangle\subset\langle V_{i}\rangle$
for all $1\leq i\leq N$ and $n\geq0$. Hence, all the assumptions of Theorem 3.8 hold for system {\rm(2.6)}
by Propositions 4.1 and 4.4. Therefore, $h(\bar{f}_{0,\infty})\geq\log\rho(A)$ by Theorem 3.8.
This completes the proof.

\bigskip

\noindent{\bf 5. Examples}\medskip

In this section, two examples are given to illustrate the theoretical results given in Sections 3
and 4.\medskip

\noindent{\bf Example 5.1.} Consider system {\rm(2.1)} with $f_n=f_1$ or $f_2$ for all $n\geq0$, where
\vspace{-0.2cm}$$f_1(x) =\begin{cases}
16x(1-x), & \; x\in [0,1/4]\cup(3/4,1], \\
                     3, &\;  x\in(1/4,3/4],\\
                     0, & \; x\in(1,3],
\end{cases}\vspace{-0.2cm}$$
and
\vspace{-0.2cm}$$f_2(x) =\begin{cases}
4x(2-x), & \; x\in [0,1/2]\cup(3/2,2], \\
                     3, &\;  x\in(1/2,3/2],\\
                     0, & \; x\in(2,3].
\end{cases}\vspace{-0.2cm}$$
Then, $\{f_n\}_{n=0}^{\infty}$ is a sequence of equi-continuous maps from $[0,3]$ to $[0,3]$.
Let
\vspace{-0.2cm}$$V_1:=[0,1/2],\;\;V_2:=[3/4,2].\vspace{-0.2cm}$$
Then, $V_1$ and $V_2$ are nonempty, disjoint, and compact subsets of $[0,3]$.
Denote
\vspace{-0.2cm}$$\mathbf{N}^{1}:=\{n\in\mathbf{N}: f_n=f_1\}.\vspace{-0.2cm}$$

In the case that $n\in\mathbf{N}^{1}$, set
\vspace{-0.2cm}$$V_{1,n}:=[0,1/4],\;\; V_{2,n}:=[3/4,1].\vspace{-0.2cm}$$
Then, $V_{1,n}\subset V_1$, $V_{2,n}\subset V_2$, and $f_{n}(V_{1,n})=f_{n}(V_{2,n})=[0,3]$. Thus,
\vspace{-0.2cm}$$
V_{1,n+1}\cup V_{2,n+1}\subset V_1\cup V_2\subset f_{n}(V_{1,n})\cap f_{n}(V_{2,n}).
\eqno(5.1)\vspace{-0.2cm}$$

In the case that $n\in\mathbf{N}\setminus\mathbf{N}^{1}$, set
\vspace{-0.2cm}$$V_{1,n}:=[0,1/2],\;\;V_{2,n}:=[3/2,2].\vspace{-0.2cm}$$
Then, $V_{1,n}\subset V_1$, $V_{2,n}\subset V_2$, and $f_{n}(V_{1,n})=f_{n}(V_{2,n})=[0,3]$. Thus,
\vspace{-0.2cm}$$
V_{1,n+1}\cup V_{2,n+1}\subset V_1\cup V_2\subset f_{n}(V_{1,n})\cap f_{n}(V_{2,n}).
\eqno(5.2)\vspace{-0.2cm}$$
It follows from  {\rm(5.1)} and {\rm(5.2)} that system {\rm(2.1)} is strictly weakly $A$-coupled expanding
in $\{V_{i,n}\}_{n=0}^{\infty}$, $i=1, 2$, with $A=(a_{ij})_{2\times 2}$ and $a_{ij}=1$, $i,j=1,2$.
It is evident that $A$ is irreducible and $\sum_{j=1}^{2}a_{i_{0}j}=2$, $i_0=1, 2$.
On the other hand, one can easily verify that
\vspace{-0.2cm}$$
d(f_n(x), f_n(y))\geq 8d(x,y),\;\forall\;x,y\in V_{i,n},\;i=1,2,\;n\in\mathbf{N}^{1},
\vspace{-0.2cm}$$
and
\vspace{-0.2cm}$$d(f_n(x), f_n(y))\geq 4d(x,y),\;\forall\;x,y\in V_{i,n},\;i=1,2,\;n\in\mathbf{N}\setminus\mathbf{N}^{1}.\vspace{-0.2cm}$$
Hence, all the assumptions of Theorems 3.7, 3.8, and 4.8 hold for system (2.1).

By Theorem 3.7, for any $n\geq0$, there exists a nonempty compact subset $\Lambda_{n}\subset\bigcup_{i=1}^{2}V_{i,n}
\subset\bigcup_{i=1}^{2}V_{i}$ with $f_{n}(\Lambda_{n})=\Lambda_{n+1}$ such that the invariant subsystem of system {\rm(2.1)} on $\{\Lambda_{n}\}_{n=0}^{\infty}$ is topologically equi-conjugate to $(\Sigma_{2}^{+}(A),\sigma_A)$, and consequently
$h(f_{0,\infty},\Lambda_0)=\log\rho(A)=\log2$. Moreover, system {\rm(2.1)} is chaotic in the strong sense of Li-Yorke,
and is also distributionally chaotic.

By Theorem 3.8, one has that $h(f_{0,\infty})\geq\log\rho(A)=\log2$.

By Theorem 4.8, for any $n\geq0$, there exists a nonempty compact subset $\mathcal{C}_{n}\subset
\bigcup_{i=1}^{2}\langle V_{i,n}\rangle$ $\subset\bigcup_{i=1}^{2}\langle V_{i}\rangle$ with $\bar{f}_{n}(\mathcal{C}_{n})
=\mathcal{C}_{n+1}$ such that the invariant subsystem of system {\rm(2.6)} on $\{\mathcal{C}_{n}\}_{n=0}^{\infty}$ is
topologically equi-conjugate to $(\Sigma_{2}^{+}(A),\sigma_A)$, and consequently $h(\bar{f}_{0,\infty},\mathcal{C}_0)
=\log\rho(A)=\log2$. Moreover, system {\rm(2.6)} is chaotic in the strong sense of Li-Yorke, and is also distributionally chaotic.\medskip

\noindent{\bf Example 5.2.} Consider the following planar non-autonomous discrete system:\medskip
\vspace{-0.2cm}$$
\left(
\begin{array}{c}
x_{1}(n+1)\\
x_{2}(n+1)\\
\end{array}
\right)
=\mathbf{F}_n\left(
\begin{array}{c}
x_{1}(n)\\
x_{2}(n)\\
\end{array}
\right), \;\;n\geq0,\eqno(5.3)\vspace{-0.2cm}$$
where $(x_1,x_2)\in\mathbf{R}^{2}$ with
\vspace{-0.2cm}$$\mathbf{F}_n\left(
\begin{array}{c}
x_{1}\\
x_{2}\\
\end{array}
\right)=\left(
\begin{array}{c}
f_{n,1}(x_1,x_2)\\
f_{n,2}(x_1,x_2)\\
\end{array}
\right),\vspace{-0.2cm}$$
\vspace{-0.2cm}$$f_{n,1}(x_{1},x_{2})=n(n+1)^{-1}\sin x_2+saw_2(12x_1),\;\;f_{n,2}(x_{1},x_{2})=f_{n,1}(x_{2},x_{1}),\vspace{-0.2cm}$$
with $saw_2(t)$ being the sawtooth function defined by
\vspace{-0.2cm}$$saw_2(t)=(-1)^{m}(t-4m),\;4m-2\leq t<4m+2,\; m\in\mathbf{Z}.\vspace{-0.2cm}$$

It is evident that $\mathbf{F}_n$ is continuous in $\mathbf{R}^{2}$ for all $n\geq0$. Set
\vspace{-0.2cm}$$V_1=[-1/6,1/6]\times[-1/6,1/6],\;V_2=[1/2,5/6]\times[1/2,5/6].\vspace{-0.2cm}$$
Clearly, $V_1$ and $V_2$ are nonempty, disjoint, and compact sets with $d(V_1,V_2)>0$.

Next, it is to show that $V_1\cup V_2\subset\mathbf{F}_n(V_1)\cap\mathbf{F}_n(V_2)$ for any fixed $n\geq0$.
For any $(x_1, x_2)\in V_1$, one has that
\vspace{-0.2cm}$$f_{n,1}(x_{1},x_{2})=n(n+1)^{-1}\sin x_2+12x_1,\;\;f_{n,2}(x_{1},x_{2})=f_{n,1}(x_{2},x_{1}).                       \eqno(5.4)\vspace{-0.2cm}$$
It follows from (5.4) that, for any $(x_1, x_2)\in V_1$ with $x_1=-1/6$,
\vspace{-0.2cm}$$f_{n,1}(x_{1},x_{2})=n(n+1)^{-1}\sin x_2-2\leq-1<-1/6,                                                          \eqno(5.5)\vspace{-0.2cm}$$
and, for any $(x_1, x_2)\in V_1$ with $x_1=1/6$,
\vspace{-0.2cm}$$f_{n,1}(x_{1},x_{2})=n(n+1)^{-1}\sin x_2+2\geq1>5/6.                                                           \eqno(5.6)\vspace{-0.2cm}$$
It then follows from (5.4) that, for any $(x_1, x_2)\in V_1$ with $x_2=-1/6$,
\vspace{-0.2cm}$$f_{n,2}(x_{1},x_{2})=n(n+1)^{-1}\sin x_1-2\leq-1<-1/6,                                                          \eqno(5.7)\vspace{-0.2cm}$$
and, for any $(x_1, x_2)\in V_1$ with $x_2=1/6$,
\vspace{-0.2cm}$$f_{n,2}(x_{1},x_{2})=n(n+1)^{-1}\sin x_1+2\geq1>5/6.                                                           \eqno(5.8)\vspace{-0.2cm}$$
Since $\mathbf{F}_n$ is continuous in $V_1$, by the intermediate value theorem and (5.5)-(5.8), one has that
$V_1\cup V_2\subset\mathbf{F}_n(V_1)$. On the other hand, for any $(x_1, x_2)\in V_2$, one has that
\vspace{-0.2cm}$$f_{n,1}(x_{1},x_{2})=n(n+1)^{-1}\sin x_2+12x_1-8,\;\;f_{n,2}(x_{1},x_{2})=f_{n,1}(x_{2},x_{1}).                     \eqno(5.9)\vspace{-0.2cm}$$
With a similar argument, one can show that $V_1\cup V_2\subset\mathbf{F}_n(V_2)$. Hence, system (5.3) is
strictly $A$-coupled-expanding in $V_1$ and $V_2$ with $A=(a_{ij})_{2\times 2}$ and $a_{ij}=1$, $i,j=1,2$.

\begin{figure}
\centering
\includegraphics[width=9cm]{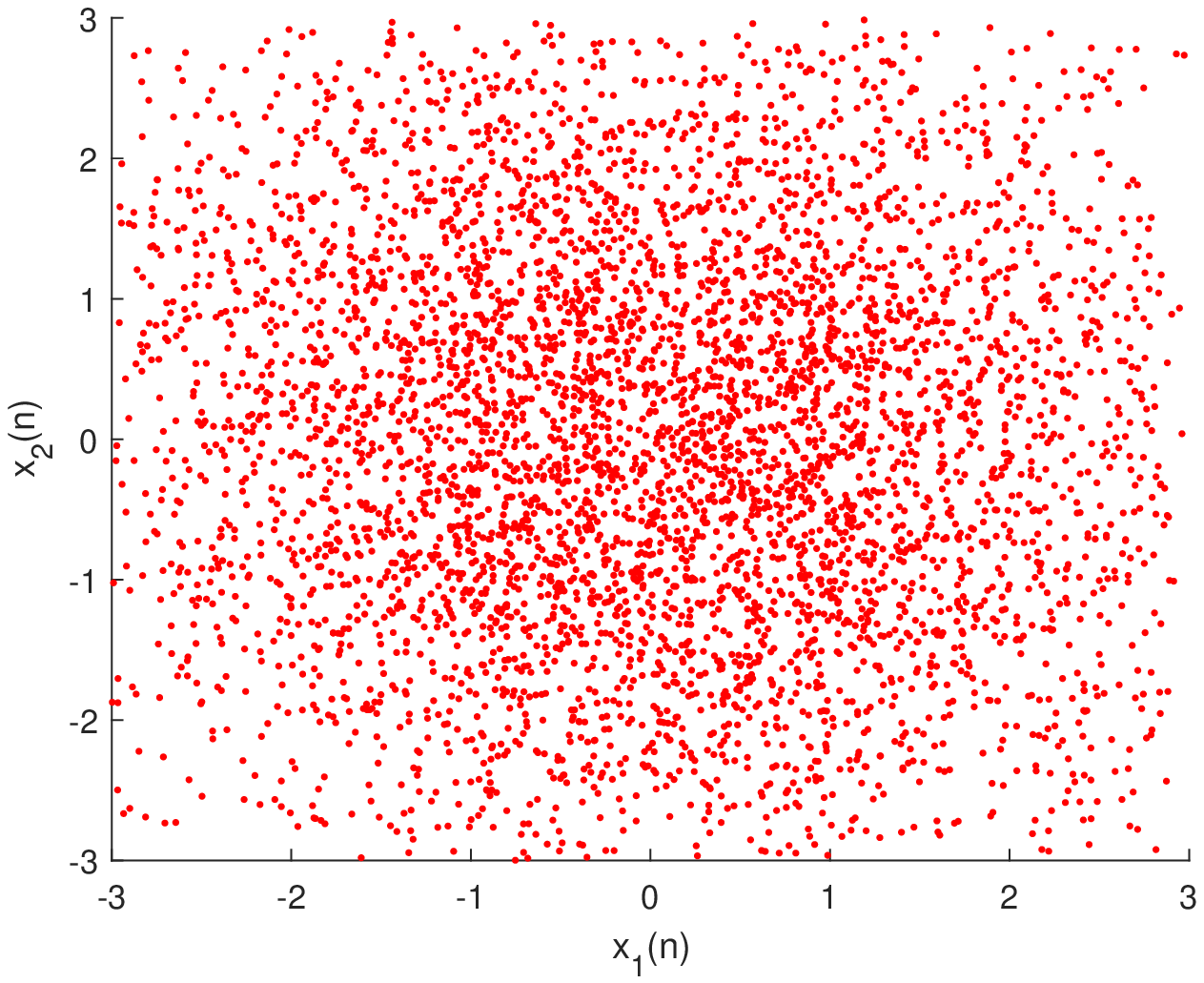}
\caption{\small{2D computer simulation result showing complex dynamical behaviors of system (5.3) in
Example 5.2 with initial conditions $x_{1}(0)=0.12,\;x_{2}(0)=0.01$, and $n=5000$.}}
%\centering
\includegraphics[width=9cm]{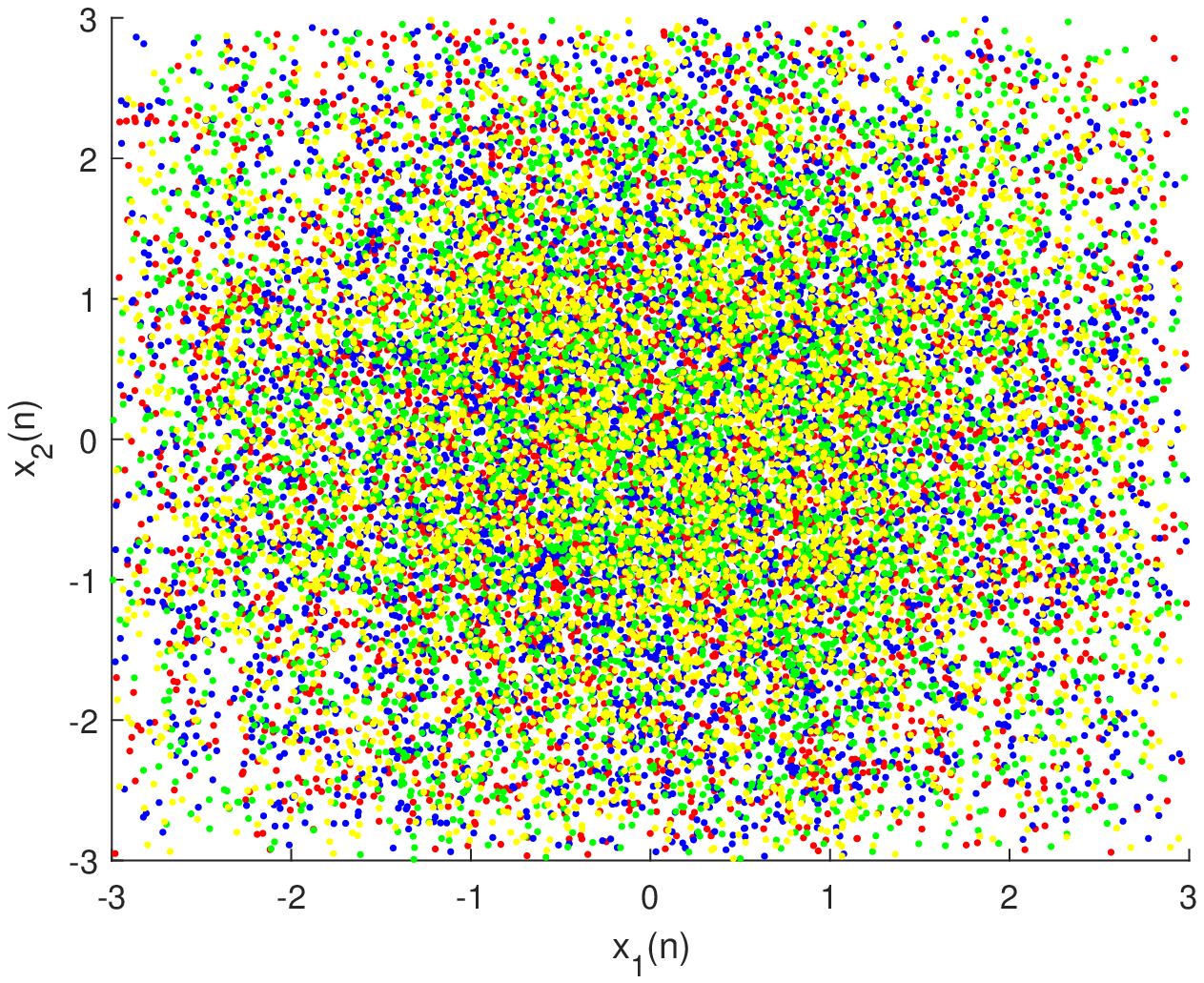}
\caption{\small{2D computer simulation result showing complex dynamical behaviors of the induced set-valued system
of system (5.3) in Example 5.2 with initial conditions
$A_0=\{(0.04,0.01),(0.05,0.01),(0.11,0.01),(0.12,0.02)\}$, and $n=5000$.}}
\end{figure}

By (5.4) and (5.9), one can easily verify that, for any $n\geq0$,
\vspace{-0.2cm}$$11\|x-y\|\leq\|\mathbf{F}_n(x)-\mathbf{F}_n(y)\|\leq13\|x-y\|,\;\forall\;x,y\in V_1\cup V_2,                      \eqno(5.10)\vspace{-0.2cm}$$
where
\vspace{-0.2cm}$$\|x\|:=\sup\{|x_i|: i=1,2\},\;\forall\;x=(x_1,x_2)\in\mathbf{R}^{2}.\vspace{-0.2cm}$$
Hence, $\{\mathbf{F}_n\}_{n=0}^{\infty}$ is equi-continuous in $\cup_{i=1}^{2}{V_i}$. Thus, all the assumptions
of Theorem 3.7 hold for system (5.3). By Theorem 3.7, one has that, for any $n\geq0$, there exists a nonempty
compact subset $\Lambda_{n}\subset\bigcup_{i=1}^{2}V_{i}$ with $\mathbf{F}_{n}(\Lambda_{n})=\Lambda_{n+1}$ such that
the invariant subsystem of system {\rm(5.3)} on $\{\Lambda_{n}\}_{n=0}^{\infty}$ is topologically equi-conjugate to
$(\Sigma_{2}^{+}(A),\sigma_A)$, and consequently system {\rm(5.3)} is chaotic in the strong sense of Li-Yorke and
is also distributionally chaotic.

Similarly, all the assumptions of Theorem 4.8 hold for the invariant subsystem of system (5.3) on
$\{\Lambda_{n}\}_{n=0}^{\infty}$. It follows from Theorem 4.8 that the induced set-valued system
of invariant subsystem of system (5.3) on $\{\Lambda_{n}\}_{n=0}^{\infty}$ has an invariant subsystem
that is topologically equi-conjugate to $(\Sigma_{2}^{+}(A),\sigma_A)$, and consequently the induced
set-valued system of system (5.3) is chaotic in the strong sense of Li-Yorke and is also distributionally chaotic.

It is remarked here that Theorem 4.8 still holds true when $f_n: X\to X$ is replaced by that
$f_n: \Lambda_n\to\Lambda_{n+1}$ with compact subsets $\Lambda_n\subset X$ for all $n\geq0$. Here,
$\bigcup_{i=1}^{2}V_{i}$ can be regarded as the base compact space $X$.\medskip

\noindent{\bf Remark 5.1.} Some simulation results are displayed in Figures 1 and 2, which show the complicated dynamical behaviors
of system (5.3) in Example 5.2 and its induced set-valued system. It is noted that the constructions of these two examples are
motivated by \cite{Huang,Shi09,Zhangli}.\medskip

\noindent{\bf Acknowledgments}\medskip

This research was supported by the Hong Kong Research Grants Council (GRF Grant CityU11200317) and the NNSF of China (Grant 11571202).
The authors would like to thank Dr. Yang Lou for his assistance in simulation.


\bigskip

\begin{thebibliography}{99}

\bibitem{birkhoff} G. D. Birkhoff, Dynamical Systems, Providence: AMS Publications, 1927.\vspace{-0.15cm}

\bibitem{block} L. Block, W. Coppel, Dynamics in One Dimension, Lecture Notes in Mathematics Vol. 1513, Springer-Verlag, Berlin/Heidelberg, 1992.\vspace{-0.15cm}

\bibitem{Devaney79} R. L. Devaney, Z. Nitecki, Shift automorphism in the H\'enon mapping, Commun. Math. Phys. 67 (1979) 137--148.\vspace{-0.15cm}

\bibitem{Devaney89} R. L. Devaney, An Introduction to Chaotic Dynamical Systems, 2nd ed, Addison-Wesley Publishing Company, 1989.\vspace{-0.15cm}

\bibitem {Fed} A. Fedeli, On chaotic set-valued discrete dynamical systems, Chaos Solit. Fract. 23 (2005) 1381--1384.\vspace{-0.15cm}

\bibitem {Gu06} R. Gu, W. Guo, On mixing property in set-valued discrete systems, Chaos Solit. Fract. 28 (2006) 747--754.\vspace{-0.15cm}

\bibitem{hadamard} J. Hadamard, Les surfaces \`{a} curbures oppos\'es et leurs lignes g\'eodesiques, J. Math. 5 (1898) 27--73.\vspace{-0.15cm}

\bibitem{Huang} Q. Huang, Y. Shi, L. Zhang, Chaotification of nonautonomous discrete dynamical systems, Int. J. Bifurcation and Chaos 21 (2011) 3359--3371.\vspace{-0.15cm}

\bibitem{ju16} H. Ju, H. Shao, Y. Choe, Y. Shi, Conditions for maps to be topologically conjugate or semi-conjugate to subshifts of finite type and
criteria of chaos, Dyn. Syst. 31 (2016) 496--505.\vspace{-0.15cm}

\bibitem{ju17} H. Ju, C. Kim, Y. Choe, M. Chen, Conditions for topologically semi-conjugacy of the induced systems to the subshift of finite type,
Chaos Solit. Fract. 98 (2017) 1--6.\vspace{-0.15cm}

\bibitem{Kennedy} J. Kennedy, J. A. Yorke, Topological horseshoes, Trans. Am. Math. Soc. 353 (2001) 2513--2530.\vspace{-0.15cm}

\bibitem{Kim} C. Kim, H. Ju, M. Chen, P. Raith, $A$-Coupled-expanding and distributional chaos, Chaos Solit. Fract. 77 (2015) 291--295.\vspace{-0.15cm}

\bibitem{kolyada} S. Kolyada, L. Snoha, Topological entropy of non-autononous dynamical systems, Random Comp. Dyn. 4 (1996) 205--233.\vspace{-0.15cm}

\bibitem{Kulczycki} M. Kulczycki, P. Oprocha, Coupled-expanding maps and matrix shifts. Int. J. Bifurcat. Chaos 23 (2013) 1--6.\vspace{-0.15cm}

\bibitem {Liao} G. Liao, X. Ma, L. Wang, Individual chaos implies collective chaos for weakly mixing discrete dynamical systems,
Chaos Solit. Fract. 32 (2007) 604--608.\vspace{-0.15cm}

\bibitem {Liu} H. Liu, E. Shi, G. Liao, Sensitivity of set-valued discrete systems, Nonlinear Anal. 71 (2009) 6122--6125.\vspace{-0.15cm}

\bibitem {Ma} X. Ma, B. Hou, G. Liao, Chaos in hyperspace system, Chaos Solit. Fract. 40 (2009) 653--660.\vspace{-0.15cm}

\bibitem {Nadler} S. B. Nadler, Continuum Theory, Pure and Applied Mathematics 158, 1992.\vspace{-0.15cm}

\bibitem{Robinson} C. Robinson, Dynamical Systems: Stability, Symbolic Dynamics and Chaos, FL: CRC Press, 1995.\vspace{-0.15cm}

\bibitem{Rom} H. Rom\'an-Flores, A note on transitivity in set-valued discrete systems, Chaos Solit. Fract. 17 (2003) 99--104.\vspace{-0.15cm}

\bibitem {san} I. S\'anchez, M. Sanchis, H. Villanueva, Chaos in hyperspaces of non-autonomous discrete systems,
Chaos Solit. Fract. 94 (2017) 68--74.\vspace{-0.15cm}

\bibitem{Shao15} H. Shao, Y. Shi, H. Zhu, Strong Li-Yorke chaos for time-varying discrete systems with A-coupled-expansion,
               Int. J. Bifurcat. Chaos 25 (2015) 1550186 (10 p.).\vspace{-0.15cm}

\bibitem{Shao16} H. Shao, Y. Shi, H. Zhu, Estimations of topological entropy for non-autonomous discrete systems,
J. Differ. Equ. Appl. 22 (2016) 474--484.\vspace{-0.15cm}

\bibitem{shaocsf} H. Shao, Y. Shi, H. Zhu, On distributional chaos in non-autonomous discrete systems, Chaos Solit. Fract.
107 (2018) 234--243.\vspace{-0.15cm}

\bibitem{shao19} H. Shao, H. Zhu, Chaos in non-autonomous discrete systems and their induced set-valued systems, Chaos 29 (2019) accepted.\vspace{-0.15cm}

\bibitem{shaocnsns} H. Shao, Y. Shi, Some weak versions of distributional chaos in non-autonomous discrete systems,
Commun. Nonlinear Sci. Numer. Simulat. 70 (2019) 318--325.\vspace{-0.15cm}

\bibitem{Shi04} Y. Shi, G. Chen, Chaos of discrete dynamical systems in complete metric spaces, Chaos Solit. Fract. 22 (2004) 555--571.\vspace{-0.15cm}

\bibitem{Shi06} Y. Shi, G. Chen, Some new criteria of chaos induced by coupled-expanding maps. In: Proc. the 1st IFAC conference on analysis
and control of chaotic systems, Reims, France, June 28-30, 2006, 157--162.\vspace{-0.15cm}

\bibitem{Shi09} Y. Shi, G. Chen, Chaos of time-varying discrete dynamical systems, J. Differ. Equ. Appl. 15 (2009) 429--449.\vspace{-0.15cm}

\bibitem{shicsf} Y. Shi, H. Ju, G. Chen, Coupled-expanding maps and one-sided symbolic dynamical systems, Chaos Solit. Fract. 39 (2009) 2138--2149.\vspace{-0.15cm}

\bibitem{smale} S. Smale, Differentiable dynamical systems, Bull. Am. Math. Soc. 73 (1967) 747--817.\vspace{-0.15cm}

\bibitem{wanghy} H. Wang, J. Xiong, Chaos for subshifts of finite type, Acta Mathematica Sinica 21 (2005) 1407--1414.\vspace{-0.15cm}

\bibitem{wang} Y. Wang, G. Wei, Conditions ensuring that hyperspace dynamical systems contain subsystems topologically (semi-)conjugate
to symbolic dynamical systems, Chaos Solit. Fract. 36 (2008) 283--289.\vspace{-0.15cm}

\bibitem{wu} X. Wu, J. Wang, G. Chen, $\mathcal{F}$-sensitivity and multi-sensitivity of hyperspatial dynamical systems,
J. Math. Anal. Appl. 429 (2015) 16--26.\vspace{-0.15cm}

\bibitem{Zhangli} L. Zhang, Y. Shi, H. Shao, Q. Huang, Chaos induced by weak A-coupled-expansion of non-autonomous discrete dynamical systems,
J. Differ. Equ. Appl. 22 (2016) 1747--1759.\vspace{-0.15cm}

\bibitem{zhangijbc} X. Zhang, Y. Shi, Coupled-expanding maps for irreducible transition matrices, Int. J. Bifurcat. Chaos 20 (2010) 3769--3783.\vspace{-0.15cm}

\bibitem{zhang} X. Zhang, Y. Shi, G. Chen, Some properties of coupled-expanding maps in compact sets, Proc. Am. Math. Soc. 141 (2013) 585--595.\vspace{-0.15cm}

\bibitem{Zhou} Z. Zhou, Symbolic Dynamics, Shanghai Scientific and Technological Education Publishing House, Shanghai, 1997.\vspace{-0.15cm}

\end{thebibliography}
\end{document}